# EVOLUTION OF DISCRETE POPULATIONS AND THE CANONICAL DIFFUSION OF ADAPTIVE DYNAMICS[1]

By Nicolas Champagnat and Amaury Lambert

*Weierstrass Institute for Applied Analysis and Stochastics and Université Pierre et Marie Curie-Paris 6 and École Normale Supérieure*

The biological theory of adaptive dynamics proposes a description of the long-term evolution of a structured asexual population. It is based on the assumptions of large population, rare mutations and small mutation steps, that lead to a deterministic ODE describing the evolution of the dominant type, called the "canonical equation of adaptive dynamics." Here, in order to include the effect of stochasticity (genetic drift), we consider self-regulated randomly fluctuating populations subject to mutation, so that the number of coexisting types may fluctuate. We apply a limit of rare mutations to these populations, while keeping the population size finite. This leads to a jump process, the so-called "trait substitution sequence," where evolution proceeds by successive invasions and fixations of mutant types. Then we apply a limit of small mutation steps (weak selection) to this jump process, that leads to a diffusion process that we call the "canonical diffusion of adaptive dynamics," in which genetic drift is combined with directional selection driven by the gradient of the fixation probability, also interpreted as an invasion fitness. Finally, we study in detail the particular case of multitype logistic branching populations and seek explicit formulae for the invasion fitness of a mutant deviating slightly from the resident type. In particular, second-order terms of the fixation probability are products of functions of the initial mutant frequency, times functions of the initial total population

Received January 2006; revised June 2006.

[1]Supported in part by a grant (Action Concertée Incitative Nouvelles Interface des Mathématiques) from the (French) Ministère de l'Enseignement Supérieur et de la Recherche, entitled "Populations Structurées" (ACI 2003-90).

*AMS 2000 subject classifications.* Primary 92D15; secondary 60K35, 60J85, 60J70, 92D10, 60J75, 92D25, 92D40, 60J80, 60J25.

*Key words and phrases.* Logistic branching process, multitype birth death competition process, population dynamics, density-dependence, competition, fixation probability, genetic drift, weak selection, adaptive dynamics, invasion fitness, timescale separation, trait substitution sequence, diffusion approximation, harmonic equations, convergence of measure-valued processes.

> size, called the invasibility coefficients of the resident by increased
> fertility, defence, aggressiveness, isolation or survival.

**1. Introduction.** Consider a multitype population where each individual of type $x$ gives birth independently at rate $b(x)$ to an individual of type $x$ (clonal reproduction) and dies either naturally at rate $d(x)$ or by *competition*. Deaths by competition are distributed as follows. Each individual of type $y$ points to each individual of type $x$ an independent exponential clock of parameter $c(x, y)$. Death by competition of an individual of type $x$ occurs as soon as a clock pointing to her rings. The vector of subpopulation sizes is distributed as the so-called *multitype logistic branching process* [23].

In this model the population size cannot go to infinity and eternal coexistence of two or more types is impossible, so the number of types decreases until it reaches 1 at some random time $T$ called the *fixation time*. For a two-type population starting with $n$ individuals of type $x$ and $m$ individuals of type $y$, we denote by $u_{n,m}(x, y)$ the probability that $y$ fixes, that is, that $y$ is the surviving type at time $T$. We refer to *neutrality* as the case when $x = y$. After time $T$ the population is said to be *monomorphic*. If the natural death rate $d$ of the surviving type is not zero, the population eventually becomes extinct, otherwise the size process is positive-recurrent.

The previous observations also apply to the population dynamics considered in this paper, that generalize the logistic branching process to a multitype setting with mutation and more general interactions, and that we call *GL*-populations. Types live in a subset $\mathcal{X}$ of $\mathbb{R}^k$ and, in a pure $x$-type population of size $n$, any individual gives birth to a new individual at rate $b(x, n)$ and dies at rate $d(x, n)$. We assume that $d(\cdot, 1) \equiv 0$ (i.e., extinction is impossible) and that, uniformly in $x$, $b(x, \cdot)$ is bounded, and $d(x, \cdot)$ bounded from below by some positive power of the total population size. Then, when there is no mutation, the size of a pure $x$-type population is positive-recurrent, as for the logistic branching process when natural death rates vanish, and converges in distribution to some random variable $\xi(x)$. However, in general, reproduction is not clonal. Specifically, each time an individual of type $x$ gives birth, the type of the newborn individual is $x$ with probability $1 - \mu(x)$ and is $x + h$ with probability $\mu(x)M(x, dh)$, where $\mu(x)$ is the *mutation probability* and $M(x, \cdot)$ is a probability measure on $\mathcal{X} - x$ called the *mutation kernel* or *mutation step law*. In this setting, the loss of diversity is counteracted by the occurrence of mutations, and the number of types can also increase.

The evolution of structured populations (with or without the presence of mutations) has long been studied, and in numerous different models. The renowned field that takes into account the complexity of the *genetic structure*—$x$ is a *genotype*—is called *population genetics* [6, 11, 15, 21].

The emerging field where the emphasis is put on the structure of *ecological interactions*—$x$ is a *phenotypic trait*, such as body size—is called *adaptive dynamics* [17, 26, 28]. The link between both is still unclear, but see [22, 31].

Let $T_n$ denote the $n$th time when the population becomes monomorphic and $V_n$ the then surviving type. The (possibly finite) sequence $(V_n)$ is of obvious interest to evolutionary biologists and is called the *trait substitution sequence*. It can be defined in much more general contexts, as soon as eternal coexistence of two or more types is not permitted by the model. It was in adaptive dynamics that this sequence was invented and studied [3, 29], under two additional assumptions. First, the biologically motivated *assumption of rare mutations* guarantees that, in the timescale of mutations (speeding up time), the widths of time intervals during which the population is polymorphic vanish, so that there is one and only one type surviving at any time $t$. To prevent the population from rapidly becoming extinct in the new timescale, one also has to rescale population sizes, thereby making the *assumption of large populations*.

Subsequently, the trait substitution sequence (TSS) is a Markov jump process on $\mathcal{X}$ whose semigroup is shown [3] to depend especially on the *invasion fitnesses* (as defined in [28]) $f(x,y)$, $x,y \in \mathcal{X}$, where $f(x,y)$ is the *expected growth rate* of a single individual of type $y$—the mutant—entering a monomorphic population of type $x$ "at equilibrium"—the residents. The evolution of the population can be described by this fitness function because of the large population assumption, which implies that deleterious mutants—those with negative invasion fitness—can never invade. Thus, evolution proceeds by *successive invasions of advantageous mutant types* replacing the resident one (selective sweeps [27]), which can be summarized by the jump process of fixed types (the TSS), also called phenotypic traits, or simply *traits*, in the adaptive dynamics setting (e.g., size, age at maturity or rate of food intake).

The TSS has been a powerful tool for understanding various evolutionary phenomena, such as evolutionary branching (evolution from a monomorphic population to a stably polymorphic one [29] that may lead to speciation [7]) and is the basis for other biological models, such as the "canonical equation of adaptive dynamics" [4, 8]. This last phrase refers to the following ODE, describing the evolution of a one-dimensional trait $x$, obtained from the TSS in the limit of *small mutations*

$$\frac{dx}{dt} = \frac{1}{2}\sigma(x)^2 \mu(x) \bar{n}(x) \frac{\partial}{\partial y} f(x,x), \tag{1}$$

where $\sigma(x)^2$ stands for the variance of the mutation step law, $\bar{n}(x)$ for the equilibrium size of a pure $x$-type population and $f(x,y)$ for the aforementioned invasion fitness. Note how only advantageous types get fixed (the trait follows the fitness gradient).

However, it is well known that slightly deleterious types can be fixed in finite populations. This phenomenon is known under the name of *genetic drift* (more generally, this term refers to allelic fluctuations which have stochastic causes). Depending on the strength of genetic drift, selection is said to be *strong* (genetic drift has negligible effects) or *weak*. Adaptive dynamics models usually assume infinite populations subject to deterministic dynamics, where only strong selection can be observed. As a consequence, the first goal of this paper is to consider *finite populations*, thereby allowing for *genetic drift* and *weak selection* [20], and continue using the bottom-up approach of adaptive dynamics; that is, model (macroscopic) evolution from (microscopic) populations [5]. In particular, we allow the population sizes to fluctuate randomly through time, in contrast to classical works in both population genetics (where the population size traditionally remains constant; see [24]) and adaptive dynamics (where population sizes are infinite).

Thus, we consider $GL$-populations in which the mutation probabilities $\mu$ are replaced with $\gamma\mu$ and time $t$ with $t/\gamma$ (mutation timescale). We prove the existence as $\gamma$ vanishes (rare mutations) of the limiting TSS and characterize its law (Theorem 3.2). Recall that $u_{n,m}(x,y)$ denotes the fixation probability of type $y$ in a two-type population starting with $n$ individuals of type $x$ and $m$ individuals of type $y$ and let $T_n$ denote the $n$th mutation time.

THEOREM 1.1. *The process $(S_t^\gamma; t \geq 0)$ defined as*

$$S_t^\gamma = \sum_{n \geq 0} V_n \mathbf{1}_{\{T_n \leq t/\gamma < T_{n+1}\}}$$

*converges in law for the Skorohod topology on $\mathbb{D}(\mathbb{R}_+, \mathcal{X})$ as $\gamma \to 0$ to a pure-jump Markov process $(S_t; t \geq 0)$, whose jumping rates $q(x, dh)$ from $x$ to $x + h$ are given by*

$$q(x, dh) = \beta(x)\chi(x, x + h)M(x, dh),$$

*where*

$$\beta(x) = \mu(x)\mathbb{E}(\xi(x)b(x, \xi(x)))$$

*and*

(2) $$\chi(x,y) = \sum_{n \geq 1} \frac{nb(x,n)\mathbb{P}(\xi(x) = n)}{\mathbb{E}(\xi(x)b(x,\xi(x)))} u_{n,1}(x,y).$$

Observe that $\beta(x)$ is the *mean mutant production rate* of a stationary $x$-type population and that $\chi(x, y)$ is the *probability of fixation* of a single (mutant) $y$-type mutant entering a pure (resident) $x$-type population with $b(x, \cdot)$-*size-biased stationary size*. In particular, $\chi(x, y)$ *is the random analogue of the usual invasion fitness*. More precisely, with the notation of (1), $\beta(x)$ is analogous to $\mu(x)b(x)\bar{n}_x$ and $\chi(x,y)$ is the exact counterpart of the expression $[f(x,y)]_+/b(y)$ obtained in [3] for the probability of fixation of a $y$-type mutant in an equilibrium $x$-type population.

The second goal of this paper is to take the limit of small mutations on the TSS defined in the last theorem. Namely, assume for simplicity that $\mathcal{X} = \mathbb{R}^k$, and denote by $Z^\epsilon$ the process $S$ modified by replacing $M(x,\cdot)$ with its image by the contraction $\epsilon$ Id and time $t$ with $t/\epsilon^2$. Under some technical assumptions which ensure that $\chi(x,\cdot)$ is differentiable, and assuming that mutation steps have zero expectation $[\int hM(x,dh) = 0]$, we get the following theorem (Theorem 4.2), where $\sigma(x)$ is the symmetric square root matrix of the covariance matrix of $M(x,\cdot)$ and $\nabla_2$ denotes the gradient taken w.r.t. the second variable.

THEOREM 1.2. *The process $Z^\epsilon$ converges in law for the Skorohod topology on $\mathbb{D}(\mathbb{R}_+, \mathbb{R}^k)$ to the diffusion process $(Z_t; t \geq 0)$ unique solution to the stochastic differential equation*

$$(3) \quad dZ_t = \beta(Z_t)\sigma^2(Z_t) \cdot \nabla_2 \chi(Z_t, Z_t)\, dt + \sqrt{\beta(Z_t)\chi(Z_t, Z_t)}\sigma(Z_t) \cdot dB_t,$$

*where $B$ is a standard $k$-dimensional Brownian motion.*

In words, we obtain a diffusion model for the evolution of a trait grounded on microscopic realistic population dynamics. The diffusion term embodies *genetic drift*. It is proportional to the mean mutant production rate $\beta(x)$, to the neutral fixation probability $\chi(x,x)$ and to the covariance matrix of $M(x,\cdot)$. The deterministic term embodies *directional selection*, and is the exact counterpart of the ODE (1).

This equation, that we christen "canonical *diffusion* of adaptive dynamics," involves the stochastic invasion fitness $\chi$ and its gradient with respect to the second variable at neutrality. As seen in the definition (2) of $\chi$, the fitness gradient only depends on the behavior of the fixation probabilities near neutrality for two interacting types (see also [24]) and on the stationary distribution $\xi(x)$. The explicit computation of these quantities is possible in the multitype logistic branching case, which we study in detail.

First consider a pure $x$-type logistic branching population with dynamical parameters $(b(x), c(x,x), d(x))$. When $d = 0$, the population size is positive-recurrent and the r.v. $\xi$ distributed according to the stationary probability is a Poisson variable of parameter $\theta := b/c$ conditioned on being nonzero. This yields

$$\beta(x) = \mu(x)b(x)\mathbb{E}(\xi(x)) = \mu(x)b(x)\frac{\theta(x)}{1 - e^{-\theta(x)}}$$

and

$$\chi(x,x) = \frac{e^{-\theta(x)} - 1 + \theta(x)}{\theta(x)^2}.$$

Second, we characterize the two-type logistic branching process by a birth vector $B$, a competition matrix $C$ and a death vector $D$:

$$B = b\mathbf{1} + \begin{pmatrix} 0 \\ \lambda \end{pmatrix},$$

$$C = c\mathbf{1} - \begin{pmatrix} 0 & 0 \\ \delta & \delta \end{pmatrix} + \begin{pmatrix} 0 & \alpha \\ 0 & \alpha \end{pmatrix} - \begin{pmatrix} 0 & \varepsilon \\ \varepsilon & 0 \end{pmatrix},$$

$$D = d\mathbf{1} - \begin{pmatrix} 0 \\ \sigma \end{pmatrix},$$

where $\mathbf{1}$ is a matrix of ones with dimensions *ad hoc* and $\lambda, \delta, \alpha, \varepsilon, \sigma$ are the selection coefficients of the mutant respectively associated to *fertility*, *defence*, *aggressiveness*, *isolation* and *survival*. In the following theorem (Theorem 5.1), we prove that each partial derivative of the fixation probability w.r.t. any selection coefficient factorizes as a function of the initial mutant frequency $p$ [either $p(1-p)$ or $p(1-p)(1-2p)$] times a function of the initial total population size, called the *invasibility coefficient* of the resident population.

THEOREM 1.3. *As a function of the multidimensional selection coefficient* $s = (\lambda, \delta, \alpha, \varepsilon, \sigma)'$, *the probability* $u$ *is differentiable, and in a neighborhood of* $s = 0$ *(selective neutrality)*,

$$u = p + v' \cdot s + o(s),$$

*where the (weak) selection gradient* $v = (v^\lambda, v^\delta, v^\alpha, v^\varepsilon, v^\sigma)'$ *can be expressed as*

$$v^\iota_{n,m} = p(1-p) g^\iota_{n+m}, \qquad \iota \neq \varepsilon,$$
$$v^\varepsilon_{n,m} = p(1-p)(1-2p) g^\varepsilon_{n+m},$$

*and the invasibility coefficients* $g$ *depend solely on the resident's characteristics* $b, c, d$ *and on the total initial population size* $n + m$.

As a result, we get

$$\nabla_2 \chi(x,x) = e^{-\theta(x)} (a_\lambda(x) \nabla b(x) - a_\delta(x) \nabla_1 c(x,x) + a_\alpha(x) \nabla_2 c(x,x)),$$

where $a_\lambda$, $a_\delta$ and $a_\alpha$ are called *adaptive slopes* in terms of respectively fertility, defence and aggressiveness. Explicit formulae are provided for the adaptive slopes as well, that can be plugged into (3). These results allow to get deeper insight into the canonical diffusion of adaptive dynamics when the microscopic interactions are those of multitype logistic branching processes.

The paper is organized as follows. First, our stochastic individual-based model is described in the next section (Section 2). Next, we state precisely the convergence results to the TSS in finite populations (Section 3) and to the canonical diffusion of adaptive dynamics (Section 4). In Section 5 we stick to the logistic branching case and obtain explicit formulae for the derivative of the fixation probability, the invasibility coefficients and the adaptive slopes. Finally (Section 6), we give the detailed proofs of the convergence results of Sections 3 and 4.

## 2. Model.

2.1. *Preliminaries.* Recall from the Introduction that a monotype (binary) logistic branching process with dynamical parameters $(b, c, d)$ is a Markov chain in continuous time $(X_t; t \geq 0)$ with nonnegative integer values and transition rates

$$q_{ij} = \begin{cases} bi, & \text{if } j = i+1, \\ ci(i-1) + di, & \text{if } j = i-1, \\ -i(b + c(i-1) + d), & \text{if } j = i, \\ 0, & \text{otherwise.} \end{cases}$$

We know from [23] that $\infty$ is an entrance boundary for $X$ and that $\mathbb{E}_\infty(\tau) < \infty$, for $\tau$ the hitting time of, say, 1. If $d \neq 0$, then the process goes extinct a.s., and if $d = 0$, it is positive-recurrent and converges in distribution to a r.v. $\xi$, where $\xi$ is a Poisson variable of parameter $\theta := b/c$ conditioned on being nonzero:

(4) $$\mathbb{P}(\xi = i) = \frac{e^{-\theta}}{1 - e^{-\theta}} \frac{\theta^i}{i!}, \qquad i \geq 1.$$

Notice that $\mathbb{E}(\xi) = \theta/(1 - \exp(-\theta))$.

For fixed $\alpha > 0$, one can generalize the interaction in the previous model to obtain the so-called $\alpha$-logistic branching process with dynamical parameters $(b, c, d)$ by modifying the transition rates as

$$q_{ij}^\alpha = \begin{cases} bi, & \text{if } j = i+1, \\ ci(i-1)^\alpha + di, & \text{if } j = i-1, \\ -i(b + c(i-1)^\alpha + d), & \text{if } j = i, \\ 0, & \text{otherwise.} \end{cases}$$

As before, it is easy to check that, when $d \neq 0$, the process goes extinct a.s., and when $d = 0$, it is positive-recurrent and converges in distribution to a r.v. $\xi^{(\alpha)}$ which law can be explicitly computed as

(5) $$\mathbb{P}(\xi^{(\alpha)} = i) = C \frac{\theta^i}{i((i-1)!)^\alpha}, \qquad i \geq 1,$$

where $\theta = b/c$ and $C$ scales the sum of these terms to 1.

2.2. *GL-populations.* In this subsection we define the general populations considered in this paper, that we call *GL*-populations, and give their basic properties. These populations are structured (multitype) populations with mutation. Their dynamics are those of eternal birth-and-death processes with birth rates at most linear and death rates of order at least $1 + \alpha$ in the total population size, hence, their name (generalized logistic).

2.2.1. *Definition.* At any time $t$, the population is composed of a finite number $N(t)$ of individuals characterized by their phenotypic traits, or simply *traits*, $x_1(t), \ldots, x_{N(t)}(t)$ belonging to a given trait space $\mathcal{X}$, assumed to be a closed subset of $\mathbb{R}^k$ for some $k \geq 1$. The population state at time $t$ will be represented by the counting measure on $\mathcal{X}$:

$$\nu_t = \sum_{i=1}^{N(t)} \delta_{x_i(t)}.$$

Let us denote by $\mathcal{M}$ the set of finite counting measures on $\mathcal{X}$, endowed with the $\sigma$-field induced by the Borel $\sigma$-field on $\mathcal{X} \subset \mathbb{R}^k$ as follows: let $\varphi$ denote the application mapping any element $\sum_{i=1}^k \delta_{x_i}$ of $\mathcal{M}$ to the $k$-tuple $(x_{\pi(1)}, \ldots, x_{\pi(k)})$ where the permutation $\pi$ of $\{1, \ldots, k\}$ is chosen such that this vector is ranked in, say, the lexicographical order. Then, this function $\varphi$ is a bijection from $\mathcal{M}$ to the set of lexicographically ordered vectors of $\bigcup_{i=0}^k \mathcal{X}^k$. The Lebesgue $\sigma$-field on this set therefore provides a $\sigma$-field on $\mathcal{M}$.

For any $\nu \in \mathcal{M}$ and any measurable function $f : \mathcal{X} \to \mathbb{R}$, we will use the notation $\langle \nu, f \rangle$ for $\int f(x) \nu(dx)$. Observe that $N(t) = \langle \nu_t, \mathbf{1} \rangle$ and that $\langle \nu_t, \mathbf{1}_\Gamma \rangle$ is the number of individuals at time $t$ with trait value in $\Gamma \subset \mathcal{X}$.

Let us consider a general structured birth-and-death process with mutation whose present state is given by the point measure $\nu$:

- $b(x, \nu)$ is the rate of birth from an individual of type $x$ in a population in state $\nu$.
- $d(x, \nu)$ is the rate of death of an individual of type $x$ in a population in state $\nu$.
- $\gamma \mu(x)$ is the probability that a birth from an individual with trait $x$ produces a mutant individual, where $\mu(x) \in [0, 1]$ and where $\gamma \in (0, 1)$ is a parameter scaling the frequency of mutations. In Section 3 we will be interested in the limit of rare mutations ($\gamma \to 0$).
- $M(x, dh)$ is the law of the trait difference $h = y - x$ between a mutant individual with trait $y$ born from an individual with trait $x$. Since the mutant trait $y = x + h$ must belong to $\mathcal{X}$, this measure has its support in $\mathcal{X} - x := \{y - x : y \in \mathcal{X}\} \subset \mathbb{R}^k$. We assume that $M(x, dh)$ has a density on $\mathbb{R}^k$ which is uniformly bounded in $x \in \mathcal{X}$ by some integrable function $\bar{M}(h)$.

In other words, the infinitesimal generator of the $\mathcal{M}$-valued Markov jump process $(\nu_t^\gamma)_{t \geq 0}$ is given by

$$
\begin{aligned}
L_\gamma \varphi(\nu) = & \int_{\mathcal{X}} [\varphi(\nu + \delta_x) - \varphi(\nu)](1 - \gamma\mu(x))b(x,\nu)\nu(dx) \\
& + \int_{\mathcal{X}} \int_{\mathbb{R}^k} [\varphi(\nu + \delta_{x+h}) - \varphi(\nu)]\gamma\mu(x)b(x,\nu)M(x,dh)\nu(dx) \\
& + \int_{\mathcal{X}} [\varphi(\nu - \delta_x) - \varphi(\nu)] \, d(x,\nu)\nu(dx).
\end{aligned}
\tag{6}
$$

We will denote by $\mathbb{P}^\gamma$ the law of this process (or $\mathbb{P}_{\nu_0}^\gamma$ when the initial condition has to be specified). When necessary, we will denote the dependence of $\nu_t$ on the parameter $\gamma$ with the notation $\nu_t^\gamma$.

DEFINITION 2.1. We say that the structured birth-and-death process with mutation defined previously is a *GL*-population, if:

- there is $\bar{b}$ such that, for any $\nu$ and $x$, $0 < b(x,\nu) \leq \bar{b}$,
- there are $\underline{c}$ and $\alpha > 0$ such that, for any $\nu$ and $x$, $\underline{c}\,(\langle \nu, \mathbf{1} \rangle - 1)^\alpha \leq d(x,\nu)$,
- if $\nu = \delta_x$, then $d(x,\nu) = 0$.

The classical structured logistic branching process [5, 9, 13, 23] is a *GL*-population with

$$b(x,\nu) = b(x)$$

and

$$d(x,\nu) = \int_{\mathcal{X}} c(x,y)(\nu - \delta_x)(dy).$$

Then, the assumption above translates as $b(\cdot) \leq \bar{b}$, $\alpha = 1$ and $c(\cdot,\cdot) \geq \underline{c}$.

2.2.2. *Basic properties.* First, the total population size $N_t = \langle \nu_t, \mathbf{1} \rangle$ of a *GL*-population is dominated by a scalar $\alpha$-logistic branching process with parameters $(\bar{b}, \underline{c}, 0)$, so $(\nu_t; t \geq 0)$ has infinite lifetime.

Second, when there is only one individual in the population ($\nu_t = \delta_x$ for some $x \in \mathcal{X}$), the death rate equals 0, so that extinction is impossible.

Third, for a *GL*-population with two types $x$ and $y$ and *no* mutation ($\mu \equiv 0$), $\nu_t = X_t \delta_x + Y_t \delta_y$, where $(X_t, Y_t : t \geq 0)$ is a bivariate Markov chain. We refer to (selective) *neutrality* as the case when $x = y$. Because of the previous domination, for the Markov chain $(X,Y)$, the union of the axes

$$\Omega_1 := \mathbb{N} \times \{0\} \quad \text{and} \quad \Omega_2 := \{0\} \times \mathbb{N}$$

is accessible and absorbing, and its complementary set is transient. So $\mathbb{P}(T < \infty) = 1$, where $T := T_{\Omega_1} \wedge T_{\Omega_2}$, and for any subset $\Gamma$ of $\mathbb{N}^2$, $T_\Gamma$ denotes the first hitting time of $\Gamma$ by $(X,Y)$. Also notice that for any $(n,m) \neq (0,0)$, $\mathbb{P}_{n,m}(T_{\Omega_1} = T_{\Omega_2}) = 0$. Then we call *fixation* (of the mutant $y$) the event $\{T_{\Omega_2} < T_{\Omega_1}\}$. The probability of fixation will be denoted by $u_{n,m}(x,y)$:

$$u_{n,m}(x,y) := \mathbb{P}(T_{\Omega_2} < T_{\Omega_1} | X_0 = n, Y_0 = m).$$

More generally, we have the following result.

PROPOSITION 2.2. *Consider a GL-population $(\nu_t; t \geq 0)$ with no mutation. For any initial condition, the fixation time $T$,*

$$T := \inf\{t \geq 0 : |\operatorname{Supp}(\nu)| = 1\},$$

*is finite a.s. and from time $T$, the population remains monomorphic with type, say, $x$. Then conditional on $T$ and $x$, the post-$T$ size process $(N(t); t \geq T)$ is positive-recurrent and converges in distribution to a random integer $\xi(x)$ such that $\sup_x \mathbb{E}(\xi(x)^n) < \infty$ for any $n$.*

PROOF. This result follows from the domination of the total population size by a $\alpha$-logistic branching process $(Z_t; t \geq 0)$ with dynamical parameters $(\bar{b}, \underline{c}, 0)$, and from the fact that such a process is positive recurrent with stationary distribution given by (5). Then, the total population size hits 1 in a.s. finite time (greater than or equal to $T$), and the population size returns then to 1 in a time bounded by the time of return to 1 of $Z$, which has finite expectation as a consequence of its positive-recurrence. So the post-$T$ size process $N$ is positive-recurrent as well and its stationary distribution is dominated by the law (5), which has finite moments. □

More generally, the domination of the population size $\langle \nu_t^\gamma, \mathbf{1} \rangle$ by a monotype $\alpha$-logistic branching process with dynamical parameters $(\bar{b}, \underline{c}, 0)$ for any $\gamma \in (0,1)$ allows us to prove the following long time bound for the moments of $\nu^\gamma$.

PROPOSITION 2.3. *Fix $p \geq 1$ and pick a positive $C$. There is a constant $C'$ such that, for any $\gamma \in (0,1)$,*

$$\mathbb{E}(\langle \nu_0^\gamma, \mathbf{1} \rangle^p) \leq C \implies \sup_{t \geq 0} \mathbb{E}(\langle \nu_t^\gamma, \mathbf{1} \rangle^p) \leq C'.$$

PROOF. With the notation of the previous proof, it suffices to show that $\sup_{t \geq 0} \mathbb{E}(Z_t^p) < +\infty$. Let us define $p_t^k = \mathbb{P}(Z_t = k)$. The backward Kolmogorov equation reads

$$\frac{d}{dt} \mathbb{E}(Z_t^p) = \sum_{k \geq 1} k^p \frac{dp_t^k}{dt}$$

$$= \sum_{k\geq 1} k^p[\bar{b}(k-1)p_t^{k-1} + \underline{c}(k+1)k^\alpha p_t^{k+1} - k(\bar{b} + \underline{c}(k-1)^\alpha)p_t^k]$$

$$= \sum_{k\geq 1}\left[\bar{b}\left(\left(1+\frac{1}{k}\right)^p - 1\right) + \underline{c}(k-1)^\alpha\left(\left(1-\frac{1}{k}\right)^p - 1\right)\right]k^{p+1}p_t^k.$$

Now, for any $k > k_0$, where $k_0 := [(2\bar{b}/\underline{c})^{1/\alpha}] + 1$, $\underline{c}(k-1)^\alpha \geq 2\bar{b}$. Therefore, for $k > k_0$,

$$\bar{b}((1+1/k)^p - 1) + \underline{c}(k-1)^\alpha((1-1/k)^p - 1)$$
$$\leq -\bar{b}[3 - 2(1-1/k)^p - (1+1/k)^p],$$

which is equivalent to $-\bar{b}p/k$. Then, enlarging $k_0$ if necessary, we obtain

$$\frac{d}{dt}\mathbb{E}(Z_t^p) \leq \sum_{k=1}^{k_0}\bar{b}(2^p-1)k_0^p - \sum_{k\geq k_0+1}\frac{\bar{b}p}{2}k^p p_t^k$$

$$\leq K - \frac{\bar{b}p}{2}\mathbb{E}(Z_t^p),$$

where the constant $K$ depends solely on $k_0$. This differential inequality yields

$$\mathbb{E}(Z_t^p) \leq \frac{2K}{\bar{b}p} + \left(\mathbb{E}(Z_0^p) - \frac{2K}{\bar{b}p}\right)e^{-\bar{b}pt/2},$$

which gives the required uniform bound. $\square$

**3. The trait substitution sequence in finite populations.** In this section we consider the $GL$-population $\nu^\gamma$ of Section 2. Our goal is to apply to this process a limit of rare mutations ($\gamma \to 0$), while keeping the population size finite, in order to describe the evolutionary process on the mutation timescale $t/\gamma$ as a "trait substitution sequence" (TSS, [3, 29]) where evolution proceeds by successive fixations of mutant types.

Let us introduce the following strong form of convergence in law. We will say that a sequence of random variables $(X_n)$ *converges strongly in law* to a r.v. $Y$ if and only if $\mathbb{E}(f(X_n)) \to \mathbb{E}(f(Y))$ as $n \to \infty$ for any bounded measurable real function $f$. Let us also make a slight abuse of notation by writing $b(x,n)$ instead of $b(x,n\delta_x)$ for the birth rate in a monomorphic population (as in the Introduction).

Fix $x \in \mathcal{X}$. For $\gamma \in (0,1)$, let the population start at $\nu_0^\gamma = N_0^\gamma \delta_x$, where the $\mathbb{N}^*$-valued random variables $N_0^\gamma$ satisfy $\sup_{\gamma \in (0,1)} \mathbb{E}((N_0^\gamma)^p) < \infty$ for some $p > 1$.

THEOREM 3.1. *For any $0 < t_1 < \cdots < t_n$, the $n$-tuple $(\nu_{t_1/\gamma}^\gamma, \ldots, \nu_{t_n/\gamma}^\gamma)$ converges strongly in law to $(\zeta_{t_1}, \ldots, \zeta_{t_n})$, where $\zeta_{t_i} = N_{t_i}\delta_{S_{t_i}}$, and*

(i) $(S_t; t \geq 0)$ *is a Markov jump process on* $\mathcal{X}$ *with initial value* $S_0 = x$ *and whose jumping rates* $q(x, dh)$ *from* $x$ *to* $x + h$ *are given by*

(7) $$q(x, dh) = \beta(x)\chi(x, x+h)M(x, dh),$$

*where*

(8) $$\beta(x) = \mu(x)\mathbb{E}(\xi(x)b(x, \xi(x)))$$

*and*

(9) $$\chi(x, y) = \sum_{n \geq 1} \frac{nb(x, n)\mathbb{P}(\xi(x) = n)}{\mathbb{E}(\xi(x)b(x, \xi(x)))} u_{n,1}(x, y).$$

(ii) *Conditional on* $(S_{t_1}, \ldots, S_{t_n}) = (x_1, \ldots, x_n)$, *the* $N_{t_i}$ *are independent and respectively distributed as* $\xi(x_i)$.

Observe that $\beta(x)$ in (8) can be seen as the *mean mutant production rate* of a stationary $x$-type population, and that $\chi(x, y)$ is the *probability of fixation* of a single $y$-type mutant entering a pure $x$-type population with $b(x, \cdot)$-*size-biased stationary size*. In particular, $\chi(x, y)$ is the random analogue of the traditional invasion fitness, in the sense proposed by Metz, Nisbet and Geritz [28].

This result shows that, in the limit of rare mutations, on the mutation timescale, the population is always monomorphic and that the dominant trait of the population evolves as a jump process over the trait space, where a jump corresponds to the appearance *and* fixation of a mutant type.

Let us denote by $\tau_n$ the $n$th mutation time ($\tau_0 = 0$), by $\rho_n$ the first time after time $\tau_n$ when the population becomes monomorphic, and by $V_n$ the single trait value surviving at time $\rho_n$ ($\rho_0 = 0$ if the initial population is monomorphic). With this notation, we can state the following result, addressing the main biological issue of Theorem 3.1, namely, the convergence of the support of the measure $\nu_{\cdot/\gamma}^\gamma$ to the process $S$.

THEOREM 3.2. *The process* $(S_t^\gamma; t \geq 0)$ *defined as*

$$S_t^\gamma = \sum_{n=0}^\infty V_n \mathbf{1}_{\{\rho_n \leq t/\gamma < \rho_{n+1}\}}$$

*converges in law for the Skorohod topology on* $\mathbb{D}(\mathbb{R}_+, \mathcal{X})$ *as* $\gamma \to 0$ *to the process* $(S_t; t \geq 0)$ *with initial state* $S_0 = x$ *characterized by* (7).

Observe that such a convergence for the measure $\nu_{\cdot/\gamma}^\gamma$ cannot hold because the population size $N_t$ in Theorem 3.1 is not a càdlàg process.

The proofs of the two preceding theorems are put to Section 6.2.

**4. The canonical diffusion of adaptive dynamics.**

4.1. *Notation and assumptions.* For any integer $r$, we denote by $\mathcal{C}_b^r$ the set of $r$ times differentiable functions with (image space *ad hoc* and) all derivatives bounded.

For a two-type *GL*-population, we will make a slight abuse of notation by writing $b(x, y, n(x), n(y))$ instead of $b(x, n(x)\delta_x + n(y)\delta_y)$ and the analogue notation for the death rate.

If $\|\cdot\|$ is the $L^1$ norm, for any $i = (n, m)$ and $j$ such that $\|j - i\| = 1$, then let $\pi_{ij}$ stand for the transition probability of the embedded Markov chain associated to the two-type *GL*-population without mutation. For example, if $j - i = (1, 0)$,

$$
\begin{aligned}
\pi_{ij}(x, y) = nb(x, y, n, m) \times (nb(x, y, n, m) + mb(y, x, m, n) \\
+ nd(x, y, n, m) + md(y, x, m, n))^{-1}.
\end{aligned}
\tag{10}
$$

From now on, we make the following additional assumptions:

- $\mathcal{X} = \mathbb{R}^k$ for simplicity,
- for all $n, m$, the functions $b(\cdot, \cdot, n, m)$ and $d(\cdot, \cdot, n, m)$ are in $\mathcal{C}_b^2$,
- there are constants $C_1, C_2$ such that, for any $x, y \in \mathcal{X}$, for any $i = (n, m)$ and $j$ such that $\|j - i\| = 1$,

$$
\|\nabla_2 \pi_{ij}(x, y)\| \leq C_1, \qquad \|H_2 \pi_{ij}(x, y)\| \leq C_2,
\tag{11}
$$

  where $\nabla_2$ is the gradient and $H_2$ the Hessian matrix taken w.r.t. the second variable,
- the mutation kernels $M(x, \cdot)$ satisfy:
  - for any $x \in \mathbb{R}^k$, $M(x, \cdot)$ has 0 expectation, that is, $\int_{\mathbb{R}^k} h M(x, dh) = 0$.
  - the covariance matrix of $M(x, \cdot)$ has Lipschitz entries and is uniformly elliptic in $x$, that is, there is a positive constant $C$ such that $\int_{\mathbb{R}^k} (s'h)^2 M(x, dh) \geq C\|s\|^2$ for any $s \in \mathbb{R}^k$.
  - the third-order moments of $M(x, \cdot)$ are uniformly bounded in $x$.

Recall that there is a symmetric matrix $\sigma(x)$ such that $\sigma(x)\sigma(x)' = \sigma(x)^2$ is the covariance matrix of $M(x, \cdot)$ which is called its square root. Then its uniform ellipticity ensures that $\sigma(x)$ has also Lipschitz entries in $x$.

4.2. *Differentiability of the probability of fixation.* In the following theorem we state the existence of the partial derivatives of $u_{n,m}(x, y)$, and show that these derivatives are always sublinear in the initial condition. We also give a uniform bound for the second-order derivatives of the fixation probability.

THEOREM 4.1. *(a) The fixation probability $y \mapsto u_{n,m}(x,y)$ is differentiable and its derivative $v_{n,m}(x,y)$ satisfies*

$$\sup_{n,m,x,y} \frac{\|v_{n,m}(x,y)\|}{n+m} < +\infty. \tag{12}$$

*(b) In addition, $y \mapsto u_{n,m}(x,y)$ is in $\mathcal{C}^2$, and its second-order derivatives are bounded by $c(n+m)^2$, where $c$ does not depend on $x$ and $y$.*

This theorem will be proved in Section 6.3.

It is easy to see from the last theorem and from the fact that $\xi(x)$ is stochastically dominated by a random variable with law (5) (with $b = \bar{b}$ and $c = \underline{c}$) that the function $\chi$ defined in (9) is in $\mathcal{C}_b^2$.

4.3. *The canonical diffusion of adaptive dynamics.* Here, we want to apply a limit of small mutation steps (weak selection) to the TSS $S$ defined in the previous section, in order to obtain the equivalent of the canonical equation of adaptive dynamics, but in finite populations [4, 8].

The limit of small jumps is obtained by introducing a parameter $\epsilon > 0$ and replacing the mutation kernels $M(x, \cdot)$ with their image by the application $h \mapsto \epsilon h$. Of course, this scaling of jumps' sizes has to be combined with a scaling of time in order to observe a nontrivial limit. This leads to the following generator for the rescaled TSS $(Z_t^\epsilon; t \geq 0)$ (as will appear further below, the factor $1/\epsilon^2$ is the right timescaling):

$$A_\epsilon \varphi(x) = \frac{1}{\epsilon^2} \int_{\mathbb{R}^k} (\varphi(x+\epsilon h) - \varphi(x))\beta(x)\chi(x, x+\epsilon h) M(x, dh). \tag{13}$$

Fix a function $\varphi$ in $\mathcal{C}_b^3$. For any $x, h \in \mathbb{R}^k$ and $\epsilon > 0$, there exists $0 \leq \epsilon_1 \leq \epsilon$ such that

$$\varphi(x+\epsilon h) - \varphi(x) = \epsilon h'\nabla\varphi(x) + \frac{\epsilon^2}{2} h' H\varphi(x+\epsilon_1 h)h,$$

where $H\varphi(y)$ denotes the Hessian matrix of $\varphi$ at $y$, and there exists $0 \leq \epsilon_2 \leq \epsilon$ such that

$$\chi(x, x+\epsilon h) = \chi(x,x) + \epsilon h' \nabla_2 \chi(x, x+\epsilon_2 h),$$

where $\nabla_2 \chi$ is the gradient of $\chi$ w.r.t. the second variable. Therefore, using the fact that $H\varphi$ and $\nabla_2 \chi$ are bounded Lipschitz functions, it takes only elementary computations to prove that

$$(\varphi(x+\epsilon h) - \varphi(x))\chi(x, x+\epsilon h)$$
$$= \epsilon(h'\nabla\varphi(x))\chi(x,x) + \epsilon^2(h'\nabla\varphi(x))(h'\nabla_2\chi(x,x))$$
$$+ \frac{\epsilon^2}{2}(h' H\varphi(x)h)\chi(x,x) + O(\epsilon^3 \|h\|^3),$$

where the $O(\epsilon^3 \|h\|^3)$ is uniform in $x \in \mathbb{R}^k$. Now, since the mutation kernel has zero expectation,

$$\int_{\mathbb{R}^k} (h'\nabla\varphi(x))\chi(x,x)M(x,dh) = 0.$$

Combining these results, thanks to boundedness of the third-order moments of the mutation kernel, we easily obtain, for any $\varphi \in \mathcal{C}_b^3$, that $A_\epsilon\varphi$ converges uniformly to the function $A_0\varphi$ defined as

$$A_0\varphi(x) = \int_{\mathbb{R}^k} (h'\nabla\varphi(x))\beta(x)(h'\nabla_2\chi(x,x))M(x,dh)$$
$$\phantom{A_0\varphi(x) =} + \tfrac{1}{2}\int_{\mathbb{R}^k} (h'H\varphi(x)h)\beta(x)\chi(x,x)M(x,dh). \tag{14}$$

In view of this, the following theorem is natural. Recall that $\sigma(x)$ is the symmetric square root matrix of the covariance matrix of $M(x,\cdot)$, which is Lipschitz in $x$.

THEOREM 4.2. *If the family $(Z_0^\epsilon)_{\epsilon>0}$ has bounded first-order moments and converges in law as $\epsilon \to 0$ to a random variable $Z_0$, then the process $Z^\epsilon$ generated by (13) with initial state $Z_0^\epsilon$ converges in law for the Skorohod topology on $\mathbb{D}(\mathbb{R}_+, \mathbb{R}^k)$ to the diffusion process $(Z_t; t \geq 0)$ with initial state $Z_0$ unique solution to the stochastic differential equation*

$$dZ_t = \beta(Z_t)\sigma^2(Z_t) \cdot \nabla_2\chi(Z_t, Z_t)\,dt + \sqrt{\beta(Z_t)\chi(Z_t, Z_t)}\sigma(Z_t) \cdot dB_t, \tag{15}$$

*where $B$ is a standard $k$-dimensional Brownian motion.*

Theorem 4.2 will be proved in Section 6.4.

Note that, by an elementary martingale (or exchangeability) argument, the neutral fixation probability $u_{n,m}(x,x)$ equals $m/(n+m)$, so that

$$\chi(x,x) = \sum_{n\geq 1} \frac{nb(x,n)\mathbb{P}(\xi(x) = n)}{\mathbb{E}(\xi(x)b(x,\xi(x)))} \frac{1}{n+1}$$
$$= \frac{1}{\mathbb{E}(\xi(x)b(x,\xi(x)))}\mathbb{E}\left(\frac{\xi(x)b(x,\xi(x))}{1+\xi(x)}\right).$$

REMARK 4.3. In the case where $\mathcal{X} \neq \mathbb{R}^k$, this result is still valid, apart from the following technical difficulties. First, for the process $Z^\epsilon$ to be well defined, one needs to assume that scaling the mutation law $M(x,dh)$ cannot drive $Z^\epsilon$ out of $\mathcal{X}$. This is true, for example, when $\epsilon \leq 1$, if $\mathcal{X}$ is convex, or if $\mathrm{Supp}(M(x,dh))$ is convex for any $x \in \mathcal{X}$. Second, uniqueness in law has to hold for the diffusion with generator $A_0$. For example, one can ensure the existence of a Lipschitz factorization $\sigma(x)\sigma(x)'$ of the covariance matrix of $M(x,\cdot)$ by assuming that the function from $\mathcal{X}$ to the set of nonnegative symmetric matrices mapping $x$ to the covariance matrix of $M(x,\cdot)$ can be extended to $\mathbb{R}^k$ in a $\mathcal{C}^2$ fashion (see [14]).

REMARK 4.4. In the case where the mutation step law $M(x,\cdot)$ has nonzero expectation, the calculation above shows that the first-order term in $A_\epsilon$ does not vanish, so that the correct timescaling is $1/\epsilon$ (instead of $1/\epsilon^2$) and the TSS process $Z^\epsilon$ can be shown to converge to the solution of the deterministic ODE

$$\frac{dz}{dt} = \beta(z)\chi(z,z) \int_{\mathbb{R}^k} hM(z,dh).$$

In this case, the main force driving evolution is the mutation bias. The mutation rate $\beta(x)$ and the fixation probability $\chi(x,x)$ only affect the speed of evolution.

Theorem 4.2 gives the equivalent of the canonical equation of adaptive dynamics [4, 8] when the population is finite. It is no longer a deterministic ODE, but a diffusion process, in which the genetic drift remains present, as a consequence of the population finiteness and of the asymptotic of weak selection ($\epsilon \to 0$).

Diffusion processes have long been used as tools for evolutionary biology, but mainly to describe the fluctuations of allelic frequencies (see, among many others, [11, 12, 19]). In sparser works [16, 25], a diffusion can describe the evolution of the dominant or mean value of a quantitative trait in a population. Our process provides such a model, but it is grounded on a microscopic precise modeling of the population dynamics, in a realistic way. In particular, the population size is not fixed and may fluctuate randomly through time but remains finite because of the density-dependence.

The diffusion part in (15) gives the strength of genetic drift: its square is proportional to the mutation rate $\beta(x)$, the neutral fixation probability $\chi(x,x)$ and the covariance matrix of the mutation step law $M(x,dh)$. As for the deterministic drift part, observe the similarity with the standard canonical equation of adaptive dynamics (1). This term gives the expression of the deterministic strength driving evolution, which is often related in macroscopic evolutionary models to a fitness gradient. In our case, the fitness is given by the function $\chi$, which appears in the deterministic part of (15) in the shape of its gradient with respect to the second variable, in a similar way as in (1). Therefore, the "hill-climbing" process of evolution occurs here, as in the classical models of adaptive dynamics, in a fitness landscape $y \mapsto \chi(x,y)$ that depends on the current state $x$ of the population.

**5. The logistic branching case.** In this section we restrict our attention to the structured logistic branching process described in the Introduction and in Section 2. Since we seek explicit expressions for the drift and diffusion parts of the canonical diffusion, we are only concerned with two-type populations with no mutation.

5.1. *Preliminaries.* Since mutations are absent, we will usually call the first type (1) resident (or wild type) and the second type (2) mutant. Then for $i, j \in \{1, 2\}$, $b_i$ (resp. $d_i$) is the birth (resp. death) rate of type $i$, and $c_{ij}$ is the competition rate felt *by* an individual of type $i$ *from* an individual of type $j$.

Formally, the two-type logistic branching process is a bivariate integer-valued continuous-time Markov process $(X_t, Y_t; t \geq 0)$ with rate matrix $Q = (q_{ij}; i \in \mathbb{N}^2, j \in \mathbb{N}^2)$, where

$$q_{ij} = \begin{cases} b_1 n, & \text{if } i = (n, m) \text{ and } j = (n+1, m), \\ b_2 m, & \text{if } i = (n, m) \text{ and } j = (n, m+1), \\ c_{11} n(n-1) + c_{12} nm + d_1 n, & \text{if } i = (n, m) \text{ and } l = (n-1, m), \\ c_{21} mn + c_{22} m(m-1) + d_2 m, & \text{if } i = (n, m) \text{ and } j = (n, m-1), \\ -r_{nm}, & \text{if } i = (n, m) \text{ and } j = (n, m), \\ 0, & \text{otherwise,} \end{cases}$$

and where the total jumping rate $r_{nm}$ is the sum of the four jumping rates from $(n, m)$.

The law of this process conditioned on fixed initial state $(n, m)$ will be denoted by $\mathbb{P}_{n,m}$. Let $B$ denote the birth vector, $C$ the competition matrix and $D$ the death vector:

$$B = \begin{pmatrix} b_1 \\ b_2 \end{pmatrix}, \qquad C = \begin{pmatrix} c_{11} & c_{12} \\ c_{21} & c_{22} \end{pmatrix}, \qquad D = \begin{pmatrix} d_1 \\ d_2 \end{pmatrix}.$$

To comply with the $GL$-population framework, we will always consider that $c_{11} c_{12} c_{21} c_{22} \neq 0$. Assuming that the presence of a mutant form does not modify the dynamical characteristics $(b, c, d)$ of the resident, we may focus on deviations from the neutral case so as to express the dynamical parameters as

$$B = b\mathbf{1} + \begin{pmatrix} 0 \\ \lambda \end{pmatrix},$$

$$C = c\mathbf{1} - \begin{pmatrix} 0 & 0 \\ \delta & \delta \end{pmatrix} + \begin{pmatrix} 0 & \alpha \\ 0 & \alpha \end{pmatrix} - \begin{pmatrix} 0 & \varepsilon \\ \varepsilon & 0 \end{pmatrix},$$

$$D = d\mathbf{1} - \begin{pmatrix} 0 \\ \sigma \end{pmatrix}.$$

In words, deviations from the neutral case are a linear combination of five *fundamental* (additive) *selection coefficients* $\lambda$, $\delta$, $\alpha$, $\varepsilon$, $\sigma$ that are chosen to be positive when they confer an advantage to the mutant. In the sequel we will see that it is indeed convenient to assess deviations to the neutral case with the help of selection coefficients in terms of:

1. *fertility* ($\lambda$, as the usual letter standing for growth rate in discrete-time deterministic models): positive $\lambda$ means increased mutant birth rate,
2. *defence* capacity ($\delta$, as in defence): positive $\delta$ means reduced competition sensitivity of mutant individuals w.r.t. the total population size,
3. *aggressiveness* ($\alpha$, as in aggressive, or attack): positive $\alpha$ means raised competition pressure exerted from any mutant individual onto the rest of the population,
4. *isolation* ($\varepsilon$, as in exclusion): positive $\varepsilon$ means lighter cross-competition between different morphs, that would lead, if harsher, to the exclusion of the less abundant one,
5. *survival* ($\sigma$, as in survival): positive $\sigma$ means reduced mutant death rate.

Under neutrality, an elementary martingale argument shows that the fixation probability equals the initial mutant frequency $p := m/(m+n)$, that is,

$$u = p.$$

The goal of the following theorem is to unveil the dependence of $u$ upon $\lambda$, $\delta$, $\alpha$, $\varepsilon$, $\sigma$ when they slightly deviate from 0. It is proved in Sections 5.3 and 5.4.

THEOREM 5.1. *As a function of the multidimensional selection coefficient $s = (\lambda, \delta, \alpha, \varepsilon, \sigma)'$, the probability $u$ is differentiable, and in a neighborhood of $s = 0$ (selective neutrality),*

(16) $$u = p + v' \cdot s + o(s),$$

*where the (weak) selection gradient $v = (v^\lambda, v^\delta, v^\alpha, v^\varepsilon, v^\sigma)'$ can be expressed as*

$$v^\iota_{n,m} = p(1-p)g^\iota_{n+m}, \qquad \iota \neq \varepsilon,$$
$$v^\varepsilon_{n,m} = p(1-p)(1-2p)g^\varepsilon_{n+m},$$

*where the $g$'s depend solely on the resident's characteristics $b, c, d$ and on the total initial population size $n+m$. They are called the invasibility coefficients.*

The invasibility coefficients of a pure resident population are interesting to study, as they provide insight as to how the fixation probability deviates from $p$ as the selection coefficients of the mutant deviate from 0. Their name is due to the fact that they only depend on the resident's characteristics and are multipliers of the mutant's selection coefficients [24]. In the next subsection we apply this result to the canonical diffusion of adaptive dynamics. The remainder of the section is then devoted to the proof of the theorem and the study of the invasibility coefficients.

5.2. *The canonical diffusion for logistic branching populations.* Consider a logistic branching population satisfying the conditions given in Section 2.2.1. Recall that $b(x)$ is the birth rate of $x$-type individuals, their death rate $d(x) \equiv 0$ and $c(x,y)$ is the competition rate felt by $x$-type individuals from $y$-type individuals. Assume as in Section 4.1 that $b(\cdot)$ and $c(\cdot, \cdot)$ are in $\mathcal{C}_b^2$. Since $c(\cdot, \cdot) \geq \underline{c}$, it is then elementary to check that (11) holds. By (4), the invasion fitness is given by

$$\chi(x,y) = \sum_{n \geq 1} \frac{n \mathbb{P}(\xi(x) = n)}{\mathbb{E}(\xi(x))} u_{n,1}(x,y) = \sum_{n \geq 1} e^{-\theta(x)} \frac{\theta(x)^{n-1}}{(n-1)!} u_{n,1}(x,y),$$

where $\theta(x) = b(x)/c(x,x)$.

Observe that

$$\left. \frac{d}{dy} \right|_{y=x} \begin{pmatrix} c(x,x) & c(x,y) \\ c(y,x) & c(y,y) \end{pmatrix} = \nabla_1 c(x,x) \begin{pmatrix} 0 & 0 \\ 1 & 1 \end{pmatrix} + \nabla_2 c(x,x) \begin{pmatrix} 0 & 1 \\ 0 & 1 \end{pmatrix},$$

where, for bivariate $f$, $\nabla_i f$ is the gradient of $f$ w.r.t. the $i$th variable ($i = 1, 2$).

Recall Theorem 4.2 and the canonical diffusion of adaptive dynamics

$$(17) \quad dZ_t = \beta(Z_t)\sigma^2(Z_t) \cdot \nabla_2 \chi(Z_t, Z_t) \, dt + \sqrt{\beta(Z_t) \chi(Z_t, Z_t)} \sigma(Z_t) \cdot dB_t,$$

where $B$ is a standard $k$-dimensional Brownian motion. From Theorem 5.1 and elementary computations, we get $\beta(x) = \mu(x)b(x)\theta(x)/(1 - e^{-\theta(x)})$ and

$$\chi(x,x) = e^{-\theta(x)} \sum_{n \geq 1} \frac{\theta(x)^{n-1}}{(n+1)(n-1)!} = \frac{e^{-\theta(x)} - 1 + \theta(x)}{\theta(x)^2}$$

and

$$\nabla_2 \chi(x,x) = e^{-\theta(x)}(a_\lambda(x) \nabla b(x) - a_\delta(x) \nabla_1 c(x,x) + a_\alpha(x) \nabla_2 c(x,x)),$$

where, for $\iota = \lambda, \delta, \alpha$,

$$(18) \qquad a_\iota(x) = \sum_{n=1}^{\infty} v_{n,1}^\iota(x) \frac{\theta(x)^{n-1}}{(n-1)!} = \sum_{n=1}^{\infty} \frac{n\theta(x)^{n-1}}{(n+1)^2(n-1)!} g_{n+1}^\iota(x),$$

and $g^\lambda$, $g^\delta$, $g^\alpha$ are the invasibility coefficients in terms of respectively fertility, defence, aggressiveness. Since the coefficients $a_\iota$ appear as factors of the gradients of the microscopic parameters $b$ and $c$ in the deterministic part of the canonical diffusion (17), we call them the *adaptive slopes*. Explicit formulae for invasibility coefficients and adaptive slopes are given in Section 5.4.

REMARK 5.2. Observe that $\varepsilon$-invasibilities do not appear in this computation, because of the symmetry between resident and mutant types in the competition kernel. One could include $\varepsilon$-invasibilities in the formula of $\nabla_2 \chi(x,x)$ by assuming a competition matrix of the form

$$\begin{pmatrix} c_1(x,x) & c_1(x,y) \\ c_2(y,x) & c_2(y,y) \end{pmatrix}$$

for some functions $c_1$ and $c_2$ coinciding on the diagonal. Such an asymmetry between resident and mutant would not be unrealistic biologically and could be explained by the resident constructing its own niche. This ecological adaptation of the resident to its medium would then result in a difference in the competition felt by $x$ from $y$ according to whether $x$ is the resident or not.

EXAMPLE. Let us consider a one-dimensional trait $x \in \mathbb{R}$ in a population undergoing symmetric competition $c(x,y) = c(|x-y|)$. This type of competition kernel has been considered in numerous earlier works; see, for example, [7]. As a consequence, $\partial c/\partial x(x,x) = \partial c/\partial y(x,x) = 0$. We may and will assume that $c(0) = 1$. We still denote by $\sigma(x)$ the standard deviation of the mutation step law $M(x,\cdot)$. Then, thanks to forthcoming Proposition 5.14 about adaptive slopes (18), the canonical diffusion of adaptive dynamics is given by

$$dZ_t = r(Z_t)\,dt + \sigma(Z_t)\mu(Z_t)^{1/2}\left(\frac{b(Z_t)}{1-e^{-b(Z_t)}} - 1\right)^{1/2} dB_t,$$

where

$$r(x) = \frac{\mu(x)\sigma(x)^2}{2}\left(1 + \frac{4}{b(x)} + \frac{b(x)-4}{1-e^{-b(x)}}\right)b'(x).$$

In forthcoming work this diffusion and other examples will be investigated.

5.3. *Fixation probability.* Here, we go back to the two-type population with no mutation and characterize the fixation probability thanks to a discrete harmonic equation (corresponding to Kolmogorov forward equations). Then we prove Theorem 5.1.

PROPOSITION 5.3. *The fixation probability $u_{n,m}$ is the unique bounded solution to*

(19)
$$\begin{aligned} (\Delta u)_{n,m} &= 0 && \text{for } (n,m) \notin \Omega_1 \cup \Omega_2, \\ u_{n,m} &= 0 && \text{for } (n,m) \in \Omega_1, \\ u_{n,m} &= 1 && \text{for } (n,m) \in \Omega_2, \end{aligned}$$

where $\Delta$ is the harmonic (its coefficients sum to zero) operator defined for any doubly indexed sequence $w$ as

$$
\begin{aligned}
(\Delta w)_{n,m} = {} & r_{n,m} w_{n,m} - b_1 n w_{n+1,m} - b_2 m w_{n,m+1} \\
& - n(c_{11}(n-1) + c_{12}m + d_1) w_{n-1,m} \\
& - m(c_{21}n + c_{22}(m-1) + d_2) w_{n,m-1}.
\end{aligned}
\tag{20}
$$

Note that in the previously displayed equation, whenever a term is not defined, the multiplying coefficient is zero. The fact that $u_{n,m}$ satisfies (19) follows from the Markov property at the first jump time of $(X, Y)$, and the uniqueness relies on Lemma 5.4 below.

LEMMA 5.4. *Consider a subset $\Gamma$ of $\mathbb{N}^2$ such that $T_\Gamma < +\infty$ $\mathbb{P}_{n,m}$-a.s. for any $n, m \geq 0$. Then, for any function $f : \Gamma \to \mathbb{R}$ such that $|f(n,m)|/(n+m+1)$ is bounded on $\Gamma$, the equation*

$$
\begin{aligned}
(\Delta h)_{n,m} &= 0 && \text{for } (n,m) \notin \Gamma, \\
h(n,m) &= f(n,m) && \text{for } (n,m) \in \Gamma
\end{aligned}
\tag{21}
$$

*admits at most one solution $h$ such that $|h(n,m)|/(n+m+1)$ is bounded.*

PROOF. It suffices to prove that (21) with $f \equiv 0$ admit $h \equiv 0$ as unique sublinear solution. Let $h$ be such a function and fix $n, m \geq 0$. Then $(h(X_{t \wedge T_\Gamma}, Y_{t \wedge T_\Gamma}); t \geq 0)$ is a $\mathbb{P}_{n,m}$-semi-martingale for $t \leq T_\Gamma$. Since, by Proposition 2.3, $\sup_{t \geq 0} \mathbb{E}_{n,m}((X_t + Y_t)^2) < +\infty$, $(h(X_t, Y_t))_{t \geq 0}$ is actually a uniformly integrable martingale. Applying the stopping theorem at time $T_\Gamma$, we get

$$0 = \mathbb{E}_{n,m}(h(X_{T_\Gamma}, Y_{T_\Gamma}) \mathbf{1}_{T_\Gamma < +\infty}) = \mathbb{E}_{n,m}(h(X_0, Y_0)) = h(n,m),$$

which completes the proof. $\square$

PROOF OF THEOREM 5.1. As seen in Proposition 5.3, the Kolmogorov forward equations translate into a discrete harmonic equation satisfied by $u$ with boundary condition 1 on $\Omega_1$, and 0 on $\Omega_2$, written as $(\Delta u)_{n,m} = 0$, where $\Delta$ is defined in (20). Combining (16) and (20), and identifying second-order terms, we get

$$
(\Delta_0 v^\iota)_{n,m} = \begin{cases} \dfrac{nm}{(n+m)(n+m+1)}, & \text{if } \iota = \lambda, \\ \dfrac{nm}{n+m}, & \text{if } \iota = \delta, \\ \dfrac{nm}{(n+m)(n+m-1)}, & \text{if } \iota = \alpha, \sigma, \\ \dfrac{nm(n-m)}{(n+m)(n+m-1)} & \text{if } \iota = \varepsilon, \end{cases}
\tag{22}
$$

where $\Delta_0$ corresponds to the neutral case of $\Delta$: for any doubly indexed $w$,

$$(\Delta_0 w)_{n,m} = (n+m)[b+c(n+m-1)+d]w_{n,m} - bn w_{n+1,m} - bm w_{n,m+1}$$
$$(23) \qquad - n[c(n+m-1)+d]w_{n-1,m} - m[c(n+m-1)+d]w_{n,m-1},$$
$$n,m \geq 0.$$

We know from Theorem 4.1 that the vector $v = (v^\lambda, v^\delta, v^\alpha, v^\varepsilon, v^\sigma)'$ is sublinear in $(n,m)$, that is, $(\|v_{n,m}\|/(n+m))_{n,m}$ is bounded. *Since the r.h.s. in (22) are all sublinear*, Lemma 5.4 ensures that $v$ *is the unique sublinear vector in $(n,m)$ solving* (22).

Thanks to this uniqueness result, it is sufficient to show that there are solutions of (22) of the following form:

$$(24) \qquad v_{n,m}^\iota = \begin{cases} \dfrac{nm}{n+m} u_{n+m}^\iota, & \text{if } \iota = \lambda, \delta, \alpha, \sigma, \\ \dfrac{nm(n-m)}{n+m} u_{n+m}^\iota, & \text{if } \iota = \varepsilon, \end{cases}$$

where, for $\iota \neq \varepsilon$, $u^\iota$ is a *bounded* real sequence indexed by $\mathbb{N} - \{0,1\}$ ($u_1^\iota$ has no effect on the values of $v_{1,0}^\iota$ and $v_{0,1}^\iota$), and $u^\varepsilon$ is a real sequence indexed by $\mathbb{N} - \{0,1,2\}$ ($u_1^\varepsilon$ and $u_2^\varepsilon$ have no effect on the corresponding values of $v_{n,m}^\varepsilon$) such that $(nu_n^\varepsilon)_n$ is bounded. The proof will then end up by writing

$$(25) \qquad g_n^\iota = \begin{cases} n u_n^\iota, & \text{if } \iota = \lambda, \delta, \alpha, \sigma,\ n \geq 2, \\ n^2 u_n^\iota, & \text{if } \iota = \varepsilon,\ n \geq 3. \end{cases}$$

In this setting (22) holds iff

$$(Lu^\lambda)_n = \frac{1}{n(n+1)} \quad \text{and} \quad (Lu^\delta)_n = \frac{1}{n}, \qquad n \geq 2,$$

$$(26) \qquad (Lu^\alpha)_n = (Lu^\sigma)_n = \frac{1}{n(n-1)}, \qquad n \geq 2,$$

$$(L'u^\varepsilon)_n = \frac{1}{n(n-1)}, \qquad n \geq 3,$$

where $L$ (resp. $L'$) is the endomorphism of the vector space $\mathcal{L}_2$ (resp. $\mathcal{L}_3$) of real sequences indexed by $\mathbb{N} - \{0,1\}$ (resp. by $\mathbb{N} - \{0,1,2\}$) defined as

$$(Lw)_n = -b\frac{n+2}{n+1}w_{n+1} + [b+c(n-1)+d]w_n$$
$$- (n-2)\left(c + \frac{d}{n-1}\right)w_{n-1}, \qquad n \geq 2,$$
$$(27)$$
$$(L'w)_n = -b\frac{n+3}{n+1}w_{n+1} + [b+c(n-1)+d]w_n$$
$$- (n-3)\left(c + \frac{d}{n-1}\right)w_{n-1}, \qquad n \geq 3.$$

The following lemma ends the proof.

LEMMA 5.5 (Existence). *There are solutions $u^\lambda, u^\delta, u^\alpha, u^\varepsilon, u^\sigma$ of (26) such that $u^\lambda, u^\delta, u^\alpha, u^\sigma$ and $(nu_n^\varepsilon)_n$ are bounded.*

This lemma will be proved in the following subsection by actually displaying explicit expressions for these solutions. □

5.4. *Invasibility coefficients and adaptive slopes.* In this subsection we give the explicit formulae for the invasibility coefficients $g^\iota$ of Theorem 5.1 and the adaptive slopes $a_\iota$ of (18).

5.4.1. *Preliminary results.* For $k \geq -2$, let $e^{(k)}$ be the sequence defined for $n \geq 2$ (3 if $k = -2$) by

$$e_n^{(k)} = \frac{1}{n+k},$$

and for $k = 2, 3$, let $\delta^{(k)}$ denote the Dirac mass at $k$:

$$\delta_n^{(k)} = \begin{cases} 1, & \text{if } n = k, \\ 0, & \text{otherwise.} \end{cases}$$

Then it is elementary to check that, for $k \geq 1$,

$$(28) \quad Le^{(k)} = -\frac{b}{k}e^{(1)} + \frac{d}{k}e^{(-1)} - b\frac{k-1}{k}e^{(k+1)}$$
$$+ (b - (k+1)c + d)e^{(k)} + (k+1)\left(c - \frac{d}{k}\right)e^{(k-1)},$$

and that

$$(29) \quad Le^{(-1)} = -2be^{(0)} + be^{(1)} + be^{(-1)} + (c+d)\delta^{(2)}.$$

Likewise, for any $k \geq 1$ and for $k = -1$,

$$(30) \quad L'e^{(k)} = -\frac{2b}{k}e^{(1)} + \frac{2d}{k}e^{(-1)} - b\frac{k-2}{k}e^{(k+1)}$$
$$+ (b - (k+1)c + d)e^{(k)} + (k+2)\left(c - \frac{d}{k}\right)e^{(k-1)},$$

and also

$$(31) \quad L'e^{(-2)} = -2(b+d)e^{(-1)} + be^{(1)} + (b+c+d)e^{(-2)} + \left(c + \frac{d}{2}\right)\delta^{(3)}.$$

Next observe that (26) can be written in the form

$$(32) \quad \begin{aligned} Lu^\lambda &= e^{(0)} - e^{(1)} \quad \text{and} \quad Lu^\delta = e^{(0)}, \\ Lu^\alpha &= Lu^\sigma = e^{(-1)} - e^{(0)}, \\ L'u^\varepsilon &= e^{(-1)} - e^{(0)}, \end{aligned}$$

so it is likely that the $u$'s can be expressed as linear combinations of the $e^{(k)}$'s. Actually, we will show they can be expressed as such linear combinations, with a potential extra additive term whose image by $L$ (resp. $L'$) is proportional to $\delta^{(2)}$ [resp. $\delta^{(3)}$]. So we end these preliminaries with displaying two sequences: one in $\mathcal{L}_2$ whose image by $L$ is $\delta^{(2)}$, and one in $\mathcal{L}_3$ whose image by $L'$ is $\delta^{(3)}$.

Assume that, at time 0, all individuals are assigned distinct labels. We denote by $\mathbf{P}_n$ the law of the logistic branching process $(b, c, d)$ starting from $n$ individuals distinctly labeled at time 0, where the label of an individual is transmitted to its offspring. In other words, under $\mathbf{P}$, we keep track of the whole descendance of each ancestral individual.

Then for $k = 2, 3$, let $T_k$ denote the first time when the total population size (i.e., the unlabeled process) is $k$. Finally, we define

$$q_n^{(k)} := \mathbf{P}_n(\text{at time } T_k, \text{ the } k \text{ living individuals have } k \text{ distinct labels}).$$

In the tree terminology, $q^{(k)}$ is the probability that all individuals in the first surviving $k$-tuple have different ancestors at time 0. In particular, $q_k^{(k)} = 1$.

LEMMA 5.6. *Let $D^{(2)} \in \mathcal{L}_2$ and $D^{(3)} \in \mathcal{L}_3$ be the sequences defined as*

$$D_n^{(2)} = \frac{q_n^{(2)}}{\kappa(n-1)}, \qquad n \geq 2,$$

$$D_n^{(3)} = \frac{q_n^{(3)}}{\kappa'(n-1)(n-2)}, \qquad n \geq 3,$$

*where*

$$\kappa = b\left(1 - \frac{2q_3^{(2)}}{3}\right) + c + d, \qquad \kappa' = \frac{b}{2}\left(1 - \frac{q_4^{(3)}}{2}\right) + c + \frac{d}{2}.$$

*Then $LD^{(2)} = \delta^{(2)}$ and $L'D^{(3)} = \delta^{(3)}$.*

Moreover, the sequences $q^{(k)}$ satisfy the following property.

LEMMA 5.7. *For any $k \geq 2$, $(q_n^{(k)})_n$ has a nonzero limit $q_\infty^{(k)}$ as $n \to \infty$.*

Proofs of Lemmas 5.6 and 5.7 are put to Section 5.5.1.

The next (sub)subsections are devoted to results related to invasibility coefficients and adaptive slopes, which prove Lemma 5.5. In Propositions 5.8, 5.9, 5.11 and 5.13, we display the solutions of (26) such that $u^\iota$ ($\iota \neq \varepsilon$) and $(nu_n^\varepsilon)_n$ are bounded (therefore proving Lemma 5.5). We also specify the behavior of each invasibility coefficient as the population size grows to infinity. Proofs of Propositions 5.11 and 5.13 are to be found in Section 5.5.2. Proposition 5.14 gives explicit expressions for the adaptive slopes.

5.4.2. *Results for the $\lambda$-invasibility.* Here, we must find a bounded sequence $u^\lambda$ in $\mathcal{L}_2$ such that $Lu^\lambda = e^{(0)} - e^{(1)}$.

Recall from Lemma 5.6 that $D^{(2)} \in \mathcal{L}_2$ is a sequence such that $LD^{(2)} = \delta^{(2)}$, and

$$D^{(2)}_n = \frac{q^{(2)}_n}{\kappa(n-1)}, \qquad n \geq 2,$$

where $q^{(2)}_n$ is the probability that the first surviving pair in the (labeled) logistic branching process $(b, c, d)$ has two distinct ancestors in the initial $n$-tuple.

Since, by (29),

$$Le^{(-1)} = -2be^{(0)} + be^{(1)} + be^{(-1)} + (c+d)\delta^{(2)},$$

and by (28),

$$Le^{(1)} = de^{(-1)} + (-2c+d)e^{(1)} + 2(c-d)e^{(0)},$$

we can readily state the following:

PROPOSITION 5.8 (Fertility). *The sequence $u^\lambda$ defined as*

(33) $$u^\lambda = -\frac{d}{2bc}e^{(-1)} + \frac{d(c+d)}{2bc}D^{(2)} + \frac{1}{2c}e^{(1)}$$

*is a bounded sequence of $\mathcal{L}_2$ such that $Lu^\lambda = e^{(0)} - e^{(1)}$. Then the invasibility coefficient $g^\lambda$ associated to fertility ($g^\lambda_n = nu^\lambda_n$) is given by*

(34) $$g^\lambda_n = -\frac{dn}{2bc(n-1)} + \frac{d(c+d)}{2bc\kappa}\frac{nq^{(2)}_n}{n-1} + \frac{n}{2c(n+1)}, \qquad n \geq 2.$$

*In particular,*

$$\lim_{n \to \infty} g^\lambda_n = \frac{b - d + d(c+d)q^{(2)}_\infty/\kappa}{2bc}.$$

5.4.3. *Results for the $\alpha$ and $\sigma$-invasibilities.* Here, we must find bounded sequences $u^\alpha$ and $u^\sigma$ in $\mathcal{L}_2$ such that $Lu^\alpha = Lu^\sigma = e^{(-1)} - e^{(0)}$. Exactly in the same way as for the $\lambda$-invasibility coefficient, we can readily make the needed statement.

PROPOSITION 5.9 (Aggressiveness, survival). *The sequences $u^\alpha$ and $u^\sigma$ defined as*

(35) $$u^\alpha = u^\sigma = \frac{2c-d}{2bc}e^{(-1)} - \frac{(2c-d)(c+d)}{2bc}D^{(2)} + \frac{1}{2c}e^{(1)}$$

are bounded sequences of $\mathcal{L}_2$ such that $Lu^\alpha = Lu^\sigma = e^{(-1)} - e^{(0)}$. Then the invasibility coefficients associated to aggressiveness $(g_n^\alpha = nu_n^\alpha)$ and survival, $(g_n^\sigma = nu_n^\sigma)$ are given by

$$(36) \quad g_n^\alpha = g_n^\sigma = \frac{(2c-d)n}{2bc(n-1)} - \frac{(2c-d)(c+d)}{2bc\kappa}\frac{nq_n^{(2)}}{n-1} + \frac{n}{2c(n+1)}, \quad n \geq 2.$$

In particular,

$$\lim_{n\to\infty} g_n^\alpha = \lim_{n\to\infty} g_n^\sigma = \frac{b + 2c - d - (2c-d)(c+d)q_\infty^{(2)}/\kappa}{2bc}.$$

5.4.4. *Results for the $\delta$-invasibility.* For $\delta$ and $\varepsilon$-invasibility coefficients, the task is mathematically more challenging. A side-effect is that we only obtain fine results in the case when the resident species has no natural death rate. This shortcoming is not very disturbing, however, because we are especially interested in precisely those populations with stationary behavior (which are those needed for applications to adaptive dynamics). From now on, we assume that $d = 0$.

Recall that we must find a bounded sequence $u^\delta$ in $\mathcal{L}_2$ such that $Lu^\delta = e^{(0)}$.

LEMMA 5.10. *Let $\Phi$ be the sequence of $\mathcal{L}_2$ defined recursively as $\Phi_2 = 1$ and*

$$(37) \quad c(n+2)\Phi_{n+1} + [b - c(n+1)]\Phi_n - b\frac{n-2}{n-1}\Phi_{n-1} = 0.$$

*Then the sequence $(n\Phi_n)_n$ converges to a nonzero finite limit $\Phi_\infty$, and the (thus well-defined) sum*

$$S := \sum_{n\geq 2} n^{-1}\Phi_n$$

*has $3c - bS = c\Phi_\infty$.*

PROPOSITION 5.11 (Defence capacity). *Define the sequence $\phi$ of $\mathcal{L}_2$ as*

$$\phi_n := \Phi_n/c\Phi_\infty \qquad n \geq 2.$$

*Then, with $\phi_1 := 1/2c$, the sequence $u^\delta$ of $\mathcal{L}_2$ defined as*

$$u^\delta := \sum_{k\geq 1} \phi_k e^{(k)}$$

*is a bounded sequence such that $Lu^\delta = e^{(0)}$. The invasibility coefficient $g^\delta$ associated to defence capacity $(g_n^\delta = nu_n^\delta)$ is given by*

$$(38) \quad g_n^\delta = \sum_{k\geq 1} \frac{n\phi_k}{n+k}, \qquad n \geq 2.$$

*In particular,*

$$g_n^\delta \sim \frac{1}{c}\ln(n) \qquad \text{as } n \to \infty.$$

The proofs of these results are given in Section 5.5.

5.4.5. *Results for the $\varepsilon$-invasibility.* Recall from Lemma 5.6 that $D^{(3)}$ is a sequence in $\mathcal{L}_3$ such that $LD^{(3)} = \delta^{(3)}$, and

$$D_n^{(3)} = \frac{q_n^{(3)}}{\kappa'(n-1)(n-2)}, \qquad n \geq 3,$$

where $q_n^{(3)}$ is the probability that the first surviving triple in the (labeled) logistic branching process $(b, c, d)$ have three distinct ancestors in the initial $n$-tuple. Now as in the previous problem ($\delta$-invasibility), we assume that the resident species has no natural death rate, that is, $d = 0$.

Here, we must find a sequence $u^\varepsilon$ in $\mathcal{L}_3$ such that $(nu_n^\varepsilon)_n$ is bounded, and $L'u^\varepsilon = e^{(-1)} - e^{(0)}$. Recall $\theta = b/c$.

LEMMA 5.12. *Let $\Psi$ be the sequence of $\mathcal{L}_3$ defined recursively as $\Psi_3 = 1$ and*

$$(39) \qquad c(n+3)\Psi_{n+1} + [b - c(n+1)]\Psi_n - b\frac{n-3}{n-1}\Psi_{n-1} = 0.$$

*Then the sequence $(n^2\Psi_n)_n$ converges to a nonzero finite limit $\Psi_\infty$, and the (thus well-defined) sums*

$$S := \sum_{n \geq 3} n^{-1}\Psi_n \quad \text{and} \quad \Sigma := \sum_{n \geq 3} \Psi_n$$

*have*

$$\Sigma + 2\theta S = \Psi_\infty + (\theta - 3)\Sigma = 5.$$

PROPOSITION 5.13 (Isolation). *Define the sequence $\psi$ of $\mathcal{L}_3$ as*

$$\psi_n := -\Psi_n/c\Psi_\infty, \qquad n \geq 3.$$

*Then, with*

$$\psi_{-2} = -\frac{1}{b(\theta+3)},$$

$$\psi_{-1} = \frac{\theta+1}{b(\theta+3)},$$

$$\psi_1 = \frac{2\theta}{3c(\theta+3)},$$

$$\psi_2 = \frac{\Sigma}{c\Psi_\infty} - \frac{2\theta+3}{3c(\theta+3)},$$

*the sequence $u^\varepsilon$ of $\mathcal{L}_3$ defined as*

$$u^\varepsilon := \sum_{k \geq -2, k \neq 0} \psi_k e^{(k)} + \frac{1}{\theta(\theta + 3)} D^{(3)}$$

*is such that $(nu_n^\varepsilon)_n$ is bounded and $L'u^\varepsilon = e^{(-1)} - e^{(0)}$. Then the invasibility coefficient $g^\varepsilon$ associated to* isolation *($g_n^\varepsilon = n^2 u_n^\varepsilon$) is given by*

$$(40) \qquad g_n^\varepsilon = \sum_{k \geq -2, k \neq 0} \frac{n^2 \psi_k}{n+k} + \frac{1}{\kappa'\theta(\theta+3)} \frac{n^2 q_n^{(3)}}{(n-1)(n-2)}, \qquad n \geq 3.$$

*In particular,*

$$g_n^\varepsilon \sim \frac{1}{c} \ln(n) \qquad \text{as } n \to \infty.$$

The proofs of these results are given in Section 5.5.

5.4.6. *Adaptive slopes.* We will denote by $q_n^{(2)}(x)$, $\kappa(x)$ and $\phi_n(x)$ the quantities appearing in Propositions 5.8, 5.9 and Proposition 5.11 when the resident population has trait $x$ [and dynamical characteristics $(b(x), c(x, x), 0)$].

The adaptive slopes appearing in the deterministic part of the canonical diffusion (15) can be expressed as follows.

PROPOSITION 5.14. *The coefficients $a_\iota$ for $\iota = \lambda, \delta, \alpha$ can be expressed in terms of the microscopic parameters $b(x), \theta(x), q_n^{(2)}(x), \kappa(x)$ and $\phi_n(x)$ as*

$$(41) \qquad a_\lambda(x) = \frac{e^{\theta(x)}(\theta(x)^2 - 3\theta(x) + 4) - \theta(x) - 4}{2b(x)\theta(x)^2},$$

$$a_\alpha(x) = \frac{e^{\theta(x)}(\theta(x)^2 - \theta(x) + 2) - \theta(x) - 2}{2b(x)\theta(x)^2}$$

$$(42) \qquad \phantom{a_\alpha(x) =} - \frac{1}{\kappa(x)\theta(x)} \sum_{n=1}^\infty \frac{n q_{n+1}^{(2)}(x)\theta(x)^{n-1}}{(n+1)!},$$

$$(43) \qquad a_\delta(x) = \sum_{k \geq 1} \frac{\phi_k(x)}{\theta(x)^{k+2}} \int_0^{\theta(x)} u^{k-1}(e^u(u^2 - u + 1) - 1)\, du$$

$$(44) \qquad \phantom{a_\delta(x)} = \frac{1}{\theta(x)^3} \int_0^{\theta(x)} (e^u(u^2 - u + 1) - 1)\pi_x\left(\frac{u}{\theta(x)}\right) du,$$

*where for $k \geq 1$,*

$$\int_0^\theta u^{k-1}(e^u(u^2 - u + 1) - 1)\, du$$

$$(45) \quad = e^\theta \left( \theta^{k+1} - (k+2)\theta^k + (k+1)^2(k-1)! \sum_{i=0}^{k-1} \frac{(-1)^i \theta^{k-i-1}}{(k-i-1)!} \right)$$

$$- (-1)^{k-1}(k+1)^2(k-1)! - \frac{\theta^k}{k}$$

and for any $v \in [0,1)$, $\pi_x(v) := \sum_{k \geq 1} \phi_k(x) v^{k-1}$. Moreover, $\pi_x$ is solution on $[0,1)$ to

$$(46) \quad u^2(1-u)\pi_x''(u) + u(\theta(x)u(1-u) + 2 - 3u)\pi_x'(u) - 2\pi_x(u) + \frac{\theta(x)}{b(x)} = 0.$$

PROOF. It follows from Proposition 5.8 that

$$a_\lambda(x) = \sum_{n=1}^\infty \frac{n\theta(x)^{n-1}}{2c(x,x)(n+2)(n+1)(n-1)!}$$

and from Proposition 5.9 that

$$a_\alpha(x) = \sum_{n=1}^\infty \frac{\theta(x)^{n-1}}{b(x)(n+1)(n-1)!}$$

$$+ \sum_{n=1}^\infty \frac{n\theta(x)^{n-1}}{2c(x,x)(n+2)(n+1)(n-1)!}$$

$$- \sum_{n=1}^\infty \frac{q_{n+1}^{(2)}(x)\theta(x)^{n-1}}{\theta(x)\kappa(x)(n+1)(n-1)!}.$$

Elementary calculations then give (41) and (42).

For the $\delta$-invasibility, using Proposition 5.11 and switching the two sums, we get

$$a_\delta(x) = \sum_{k \geq 1} \phi_k(x) \sum_{n \geq 1} \frac{n^2 \theta(x)^{n-1}}{(n+k+1)(n+1)!}.$$

The following observation

$$\sum_{n \geq 1} \frac{n^2 u^{n-1}}{(n+1)!} = \frac{d}{du}\left( u \frac{d}{du} \sum_{n \geq 1} \frac{u^n}{(n+1)!} \right)$$

$$= \frac{d}{du}\left( u \frac{d}{du} \left( \frac{e^u - 1}{u} \right) \right)$$

$$= \frac{e^u(u^2 - u + 1) - 1}{u^2}$$

yields (43). Equation (44) follows from switching the sum and the integral in (43), which is standard since $\phi_n(x) = O(1/n)$.

Equation (45) can be checked using the fact that $e^u(u^k - ku^{k-1} + k(k-1)u^{k-2} + \cdots + (-1)^k k!)$ is a primitive of $e^u u^k$.

Finally, (46) can be deduced from the facts that $\phi_1(x) = 1/2c(x,x) = \theta(x)/2b(x)$ and

$$\forall n \geq 2 \quad (n+2)\phi_{n+1}(x) + (\theta(x) - n - 1)\phi_n(x) - \theta(x)\frac{n-2}{n-1}\phi_{n-1}(x) = 0.$$

Multiplying these equations by $(n-1)x^{n+1}$ and summing over $n \geq 2$ yields

$$\begin{aligned}
0 &= \sum_{n \geq 3}(n+1)(n-2)\phi_n(x)u^n + \theta(x)\sum_{n \geq 2}(n-1)\phi_n(x)u^{n+1} \\
&\quad - \sum_{n \geq 2}(n+1)(n-1)\phi_n(x)u^{n+1} - \theta(x)\sum_{n \geq 2}(n-1)\phi_n(x)u^{n+2} \\
&= \frac{d}{du}\left(u^4 \frac{d}{du}\left(\frac{\pi_x(u) - \phi_1(x) - \phi_2(x)u}{u}\right)\right) \\
&\quad + \theta(x)u^3\pi'_x(u) - u\frac{d}{du}(u^3\pi'_x(u)) - \theta(x)u^4\pi'_x(u),
\end{aligned}$$

which finally gives (46). $\square$

### 5.5. Proofs.

#### 5.5.1. Proofs of Lemmas 5.6 and 5.7.

PROOF OF LEMMA 5.6. It is quite elementary to check the result by standard applications of the Markov property under $\mathbf{P}$, but we prefer to give a more conceptual proof. We start with $D^{(2)}$.

Under $\mathbb{P}$, we only keep track of two types at time $t$, that is, the number $X_t$ of residents, and the number $Y_t$ of mutants, whereas under $\mathbf{P}_n$, there are $n$ types at time 0, say, $1, 2, \ldots, n$. Recall individuals of all types are exchangeable (because in this setting, the discrete operators $\Delta_0$ and $L$ are associated to selective neutrality). Set

$$w_{n,m} := \mathbb{P}_{n,m}(X_{T_2} = Y_{T_2} = 1).$$

Now, by exchangeability,

$$w_{n,m} = \sum_{i=1}^{n}\sum_{j=n+1}^{n+m} \mathbf{P}_{n+m}(\text{at } T_2, \text{ type } i \text{ and type } j \text{ have one representative each})$$
$$= nm\mathbf{P}_{n+m}(\text{at } T_2, \text{ type 1 and type 2 have one representative each}),$$

and once again, by exchangeability,

$$q_n^{(2)} = \sum_{1 \leq i < j \leq n} \mathbf{P}_n(\text{at } T_2, \text{ type } i \text{ and type } j \text{ have one representative each})$$

$$= \binom{n}{2} \mathbf{P}_n(\text{at } T_2, \text{ type 1 and type 2 have one representative each}).$$

As a consequence,

$$w_{n,m} = \frac{2nm}{(m+n)(m+n-1)} q_{m+n}^{(2)}.$$

Observe that by definition $w$ is harmonic (in the sense that $\Delta_0 w = 0$) on the complementary of $\Omega_1 \cup \Omega_2 \cup \{(1,1)\}$. Then as in the previous subsection, with $v_n = q_n^{(2)}/(n-1)$, we get that $(Lv)_n = 0$ for any $n \geq 3$. The proof is completed by checking that $(Lv)_2 = \kappa \neq 0$.

As to $D^{(3)}$, set

$$w_{n,m} := \mathbb{E}_{n,m}(X_{T_3} Y_{T_3}(X_{T_3} - Y_{T_3}))$$
$$= 2\mathbb{P}_{n,m}((X_{T_3}, Y_{T_3}) = (2,1)) - 2\mathbb{P}_{n,m}((X_{T_3}, Y_{T_3}) = (1,2)).$$

Now, by exchangeability,

$$\mathbb{P}_{n,m}((X_{T_3}, Y_{T_3}) = (1,2))$$

$$= \sum_{i=1}^{n} \sum_{j=n+1}^{n+m} \mathbf{P}_{n+m}(\text{at } T_3, \text{ type } i \text{ has one representative}$$
$$\text{and type } j \text{ has two})$$

$$+ \sum_{1 \leq i \leq n} \sum_{n+1 \leq j < k \leq n+m} \mathbf{P}_{n+m}(\text{at } T_3, \text{ types } i, j \text{ and } k \text{ have}$$
$$\text{one representative each})$$

$$= nm \mathbf{P}_{n+m}(\text{at } T_3, \text{ type 1 has one}$$
$$\text{representative and type 2 has two})$$

$$+ n\frac{m(m-1)}{2} \mathbf{P}_{n+m}(\text{at } T_3, \text{ types 1, 2 and 3 have one}$$
$$\text{representative each}).$$

Since $\mathbb{P}_{n,m}((X_{T_3}, Y_{T_3}) = (1,2)) = \mathbb{P}_{m,n}((X_{T_3}, Y_{T_3}) = (2,1))$, the corresponding first terms in the difference cancel out, and we are left with

$$w_{n,m} = nm(n-m)\mathbf{P}_{n+m}(\text{at } T_3, \text{ types 1, 2 and 3 have one representative each}).$$

But again we get an expression involving the last displayed probability as

$$q_n^{(3)} = \binom{n}{3} \mathbf{P}_n(\text{at } T_3, \text{ types 1, 2 and 3 have one representative each}),$$

so that

$$w_{n,m} = \frac{6nm(n-m)}{(m+n)(m+n-1)(m+n-2)} q^{(3)}_{m+n}.$$

This time $w$ is harmonic on the complementary of $\Omega_1 \cup \Omega_2 \cup \{(1,1),(1,2),(2,1)\}$. Then with $v_n = q^{(3)}_n/(n-1)(n-2)$, we get that $(L'v)_n = 0$ for any $n \geq 4$. The proof is completed by checking that $(L'v)_3 = \kappa' \neq 0$. □

PROOF OF LEMMA 5.7. From [23], we know that $\infty$ is an entrance boundary for the probabilities $\mathbf{P}_n$, $n \geq 1$, so that $\mathbf{P}_\infty$ and $q^{(k)}_\infty$ are properly defined for any $k \geq 1$. At time $t$, we denote by $Z_t$ the number of living individuals and by $N_t$ the number of types represented. Obviously, under $\mathbf{P}_\infty$, $Z_t \to \infty$ as $t \to 0+$. As to $N$, since it is a nonincreasing function of time, it has a right-limit $N_{0+} \leq \infty$ at $t = 0$. Next we want to show that for any $k_0 \geq 2$,

$$(47) \qquad q^{(k_0)}_\infty = 0 \implies \mathbf{P}_\infty(N_{0+} \leq k_0) = 1.$$

This will end the proof of the lemma. Indeed, $N_{0+} \leq k_0$ means that, under $\mathbf{P}_\infty$, there are at most $k_0$ individuals whose total descendance at any time $t$ is $Z_t$. Then, conditional on these individuals, $Z$ would be dominated by a binary logistic branching process starting at $k_0$, which contradicts the fact that $Z_{0+} = +\infty$. Conclude by summing over all possible $k_0$-tuples.

Now, we prove (47). Assume there is $k_0 \geq 2$ such that $q^{(k_0)}_\infty = 0$. Since for $k \geq k_0$, $q^{(k_0)}_n > q^{(k)}_n q^{(k_0)}_k$, we get that $q^{(k)}_\infty = 0$ for all $k \geq k_0$. Recall that $T_j$ is the first hitting time of $j$ by $Z$. For $n \geq j \geq k \geq k_0$,

$$q^{(k)}_n > \mathbf{P}_n(N_{T_j} = k, N_{T_k} = k)$$
$$> \mathbf{P}_n(N_{T_j} = k)C(j,k),$$

where $C(j,k)$ is the probability that, conditional on $T_j = k$, after picking $k$ representative individuals at $T_j$ (one for each type), the first $j - k$ events after $T_j$ are the deaths of all nonrepresentative individuals. Because this probability only depends on $j$ and $k$, we get that $\mathbf{P}_\infty(N_{T_j} = k) = 0$ for all $j \geq k \geq k_0$. As a consequence,

$$\mathbf{P}_\infty(N_{T_j} \leq k_0) = 1.$$

But under $\mathbf{P}_\infty$, $\lim_{j \to \infty} T_j = 0$ a.s., so that $\mathbf{P}_\infty(N_{0+} \leq k_0) = 1$. □

5.5.2. *Proofs for $\delta$ and $\varepsilon$-invasibilities.*

PROOF OF LEMMA 5.10. The proof stems immediately from the two following claims. Claim 1 will also be helpful in the proof of Proposition 5.11. □

CLAIM 1. *Assume that $(n\Phi_n)_n$ converges to a finite limit $\Phi_\infty$, where $(\Phi_n)_n$ is defined in (37). Then the real number $S := \sum_{n \geq 2} n^{-1}\Phi_n$ and the sequence $W := \sum_{k \geq 2} \Phi_k e^{(k)}$ of $\mathcal{L}_2$ are well defined, and*

(i) $c\Phi_\infty = 3c - bS$,
(ii) $LW = c\Phi_\infty e^{(1)}$.

CLAIM 2. *The sequence $(n\Phi_n)_n$ converges to a nonzero finite limit.*

PROOF OF CLAIM 1. To prove (i), let

(48) $$\beta_n := (n+1)\Phi_n \quad \text{and} \quad \gamma_n := (n-1)\Phi_n, \quad n \geq 2,$$

so that

$$\lim_{n \to \infty} \beta_n = \lim_{n \to \infty} \gamma_n = \Phi_\infty,$$

and, thanks to (37),

(49) $$\beta_{n+1} - \beta_n = -\frac{\theta}{n-1}(\gamma_n - \gamma_{n-1}), \quad n \geq 2,$$

with $\gamma_1 = 0$. As a consequence, by Abel's transform, we get

$$S = \sum_{n \geq 2} \gamma_n \left(\frac{1}{n-1} - \frac{1}{n}\right) = \sum_{n \geq 2} \frac{\gamma_n - \gamma_{n-1}}{n-1}$$

$$= -\theta^{-1} \sum_{n \geq 2} (\beta_{n+1} - \beta_n) = -\theta^{-1}(\Phi_\infty - \beta_2) = (3 - \Phi_\infty)/\theta.$$

As for (ii), thanks to (28) and (29), and by continuity of linear operators,

$$LW = \lim_{l \to \infty} \sum_{k=2}^{l} \Phi_k L e^{(k)}$$

$$= \lim_{l \to \infty} \sum_{k=2}^{l} \Phi_k \left(-\frac{b}{k}e^{(1)} - b\frac{k-1}{k}e^{(k+1)}\right.$$

$$\left. + [b - (k+1)c]e^{(k)} + (k+1)ce^{(k-1)}\right)$$

$$= \lim_{l \to \infty} \left\{-b\left(\sum_{k=2}^{l} k^{-1}\Phi_k\right)e^{(1)} - b\sum_{i=2}^{l+1}\frac{i-2}{i-1}\Phi_{i-1}e^{(i)}\right.$$

$$\left. + \sum_{k=2}^{l}[b - (k+1)c]\Phi_k e^{(k)} + \sum_{j=1}^{l-1} c(j+2)\Phi_{j+1}e^{(j)}\right\}$$

$$= -bSe^{(1)}$$

$$+ \lim_{l \to \infty} \left\{ \sum_{k=2}^{l-1} \left( -b\frac{k-2}{k-1}\Phi_{k-1} + [b - (k+1)c]\Phi_k + c(k+2)\Phi_{k+1} \right) e^{(k)} \right.$$
$$- b\frac{l-2}{l-1}\Phi_{l-1}e^{(l)}$$
$$\left. - b\frac{l-1}{l}\Phi_l e^{(l+1)} + [b - (l+1)c]\Phi_l + 3c\Phi_2 e^{(1)} \right\}$$
$$= (3c - bS)e^{(1)},$$

which ends the proof. □

PROOF OF CLAIM 2. We split this proof into the four following steps [recall (48)]:

(i) if $(n\Phi_n)_n$ converges to a finite limit $\Phi_\infty$, then $\Phi_\infty \neq 0$,
(ii) $(\beta_n)_n$ has constant sign for large $n$,
(iii) $(\beta_n)_n$ is bounded,
(iv) $(\beta_n)_n$ converges.

Since we are only interested in the asymptotic properties of the sequences $(\Phi_n)_n$, $(\beta_n)_n$ and $(\gamma_n)_n$, we will implicitly assume throughout this proof that $\theta/(n+1) < 1$, that is, $n \geq \theta$.

(i) If $\Phi_\infty$ exists, then thanks to Claim 1, we can define $W = \sum_{k \geq 2} \Phi_k e^{(k)}$ and the doubly indexed sequence $w$ as

$$w_{n,m} = \frac{nm}{n+m}W_{n+m}, \qquad (n,m) \in \mathbb{N} \times \mathbb{N} \setminus (0,0).$$

Because $W$ is bounded, $w$ is sublinear. Assume $\Phi_\infty = 0$. Then by Claim 1, $LW = 0$, and the same calculations as those yielding (26) and (27) show that $\Delta_0 w = 0$. The contradiction comes with Lemma 5.4, which implies that the null sequence is the only sublinear doubly indexed sequence which vanishes on $\Omega_1 \cup \Omega_2$ and is in the kernel of $\Delta_0$.

(ii) First observe that (37) reads

$$(50) \qquad \beta_{n+1} = \left(1 - \frac{\theta}{n+1}\right)\beta_n + \theta\frac{n-2}{n(n-1)}\beta_{n-1},$$

so if there is $n_0 \geq \theta$ such that $\beta_{n_0}\beta_{n_0-1} \geq 0$, then a straightforward induction shows that $(\beta_n)_{n \geq n_0}$ never changes sign. Now, we prove that if no such $n_0$ exists, then $(\beta_n)_n$ converges to 0, which contradicts (i). Indeed, assume that for all $n \geq \theta$, $\beta_n\beta_{n-1} < 0$. Then for any $n \geq \theta + 1$, if $\beta_{n-1} < 0$, then $\beta_n > 0$ and $\beta_{n+1} < 0$, so that

$$\left(1 - \frac{\theta}{n+1}\right)\beta_n < -\theta\frac{n-2}{n(n-1)}\beta_{n-1},$$

which can be written as
$$\left(1 - \frac{\theta}{n+1}\right)|\beta_n| < \theta\frac{n-2}{n(n-1)}|\beta_{n-1}|,$$
and we would get the same inequality if $\beta_{n-1} > 0$. This would imply that $|\beta_n|/|\beta_{n-1}|$ would vanish as $n$ grows, and so would $\beta_n$.

(iii) Without loss of generality, we can assume thanks to (ii) that there is $n_0$ such that $\beta_n \geq 0$ for all $n \geq n_0$ (otherwise change $\beta$ for $-\beta$). Next, we prove that, for all $n \geq n_0$, $\beta_{n+1} < \max(\beta_n, \beta_{n-1})$. It is then elementary to see that $(\beta_n)_n$ is bounded. First check that
$$\beta_{n+1} - \beta_n = -\frac{\theta}{n-1}(\beta_n - \beta_{n-1}) - 2\theta\frac{\beta_{n-1}}{n(n-1)(n+1)},$$
so for any $n \geq n_0$,
$$\beta_{n+1} - \beta_n \leq -\frac{\theta}{n-1}(\beta_n - \beta_{n-1}).$$
In particular, if $\beta_{n+1} \geq \beta_n$, then $\beta_n \leq \beta_{n-1}$, and
$$|\beta_{n+1} - \beta_n| \leq \frac{\theta}{n-1}|\beta_n - \beta_{n-1}| < |\beta_n - \beta_{n-1}|,$$
which reads $\beta_{n+1} - \beta_n < -\beta_n + \beta_{n-1}$, that is, $\beta_{n+1} < \beta_{n-1}$. As a conclusion, $\beta_{n+1} < \beta_n$ or $\beta_{n+1} < \beta_{n-1}$.

(iv) By (49) and Abel's transform, we get
$$\beta_{n+1} - \beta_2 = -\theta \sum_{k=2}^{n} \frac{\gamma_k - \gamma_{k-1}}{k-1}$$
$$= -\theta\frac{\gamma_n}{n-1} - \theta \sum_{k=2}^{n-1} \frac{\gamma_k}{k(k-1)},$$
and the r.h.s. converges, because $(\gamma_n)_n$ is bounded, thanks to (iii). □

PROOF OF PROPOSITION 5.11. Thanks to Claim 1 above, since $u^\delta = \phi_1 e^{(1)} + (c\Phi_\infty)^{-1}W$,
$$Lu^\delta = \phi_1 Le^{(1)} + (c\Phi_\infty)^{-1}LW$$
$$= (2c)^{-1}(-2ce^{(1)} + 2ce^{(0)}) + (c\Phi_\infty)^{-1}c\Phi_\infty e^{(1)}$$
$$= e^{(0)}.$$
The boundedness of $u^\delta$ is straightforward. To get the equivalent of $g_n^\delta$ as $n \to \infty$, it is sufficient to prove that $nW_n \sim \Phi_\infty \ln(n)$. First, starting over from the proof of Claim 2(iv) above, get that
$$\Phi_\infty - \beta_2 = -\theta \sum_{k \geq 2} \frac{\gamma_k}{k(k-1)},$$

so that

$$\beta_{n+1} - \Phi_\infty = -\theta \frac{\gamma_n}{n-1} + \theta \sum_{k \geq n} \frac{\gamma_k}{k(k-1)},$$

which implies that $\beta_n - \Phi_\infty = o(n^{-1})$. Next, writing $\rho_k := k\Phi_k$, we get

$$nW_n = \sum_{k \geq 2} \frac{n\rho_k}{k(n+k)} = \sum_{k \geq 2} \rho_k \left(\frac{1}{k} - \frac{1}{n+k}\right)$$

$$= \lim_{l \to \infty} \left\{ \sum_{k=2}^{n+1} \frac{\rho_k}{k} + \sum_{k=2}^{l} \frac{\rho_{n+k} - \rho_k}{n+k} - \sum_{k=l+1}^{l+n} \frac{\rho_k}{k} \right\}$$

$$= \sum_{k=2}^{n+1} \frac{\rho_k - \Phi_\infty}{k} + \sum_{k=2}^{n+1} \frac{\Phi_\infty}{k} + \sum_{k \geq 2} \frac{\rho_{n+k} - \rho_k}{n+k}$$

$$= \Phi_\infty \ln(n) + O(1),$$

where the last equation comes from the fact that $\rho_k = \Phi_\infty + O(k^{-1})$ as $k \to \infty$. □

PROOF OF LEMMA 5.12. Since proofs for the isolation $\varepsilon$ are quite similar to those done for the defence capacity $\delta$, we will often sketch them. The proof of Lemma 5.12 stems immediately from the following two claims. □

CLAIM 1. *Assume that $(n^2 \Psi_n)_n$ converges to a finite limit $\Psi_\infty$, where $(\Psi_n)_n$ is defined in (39). Then the real numbers $S := \sum_{n \geq 3} n^{-1} \Psi_n$ and $\Sigma := \sum_{n \geq 3} \Psi_n$, as well as the sequence $Z := \sum_{k \geq 3} \Psi_k e^{(k)}$ of $\mathcal{L}_3$ are well defined, and:*

(i) $\Sigma + 2\theta S = \Psi_\infty + (\theta - 3)\Sigma = 5.$
(ii) $L'(Z - \Sigma e^{(2)}) = c\Psi_\infty(e^{(2)} - e^{(1)}).$

CLAIM 2. *The sequence $(n^2 \Psi_n)_n$ converges to a nonzero finite limit.*

PROOF OF CLAIM 1. To prove (i), let

(51) $\beta_n := (n+2)(n+1)\Psi_n$ and $\gamma_n := (n-2)(n-1)\Psi_n, \quad n \geq 3,$

so that

$$\lim_{n \to \infty} \beta_n = \lim_{n \to \infty} \gamma_n = \Psi_\infty,$$

and, thanks to (39),

(52) $\qquad \beta_{n+1} - \beta_n = -\frac{\theta(n+2)}{(n-1)(n-2)}(\gamma_n - \gamma_{n-1}), \quad n \geq 3,$

with $\gamma_2 = 0$. As a consequence, by two applications of Abel's transform, we get

$$\Sigma = \sum_{n \geq 3} \Psi_n = \sum_{n \geq 3} \beta_n \left( \frac{1}{n+1} - \frac{1}{n+2} \right)$$

$$= \frac{\beta_3}{4} + \sum_{n \geq 3} \frac{\beta_{n+1} - \beta_n}{n+2} = 5\Psi_3 - \theta \sum_{n \geq 3} \frac{\gamma_n - \gamma_{n-1}}{(n-1)(n-2)}$$

$$= 5 - \theta \sum_{n \geq 3} \frac{2\gamma_n}{n(n-1)(n-2)}$$

$$= 5 - 2\theta S.$$

On the other hand, the same type of arguments as above show that

$$\Sigma = \sum_{n \geq 3} \frac{\gamma_n}{(n-2)(n-1)} = \sum_{n \geq 3} \frac{\gamma_n - \gamma_{n-1}}{n-2}$$

$$= -\theta^{-1} \sum_{n \geq 3} \frac{n-1}{n+2} (\beta_{n+1} - \beta_n)$$

$$= -\theta^{-1} \left\{ -\frac{\beta_3}{4} + \Psi_\infty + \sum_{n \geq 3} \left( \frac{n-2}{n+1} \beta_n - \frac{n-1}{n} \beta_n \right) \right\}$$

$$= -\theta^{-1}(-5 + \Psi_\infty - 3\Sigma),$$

which ends the proof of (i). With the help of (30) and (31), (ii) can be proved easily mimicking what was done for the $\delta$-invasibility. $\square$

PROOF OF CLAIM 2. We proceed just as for the defence capacity. First, we prove that if $\Psi_\infty$ exists, it cannot be 0. Indeed, consider $\Sigma$ and $Z \in \mathcal{L}_3$ defined in Claim 1, and further define

$$z_{n,m} = \frac{nm(n-m)}{n+m} \left( Z_{n+m} - \frac{\Sigma}{n+m+2} \right), \qquad (n,m) \in \mathbb{N}^\star \times \mathbb{N}^\star.$$

Because $(nZ_n)_n$ is bounded, $z$ is sublinear. If $\Psi_\infty = 0$, then thanks to Claim 1, we would get $\Delta_0 z = 0$, but this would contradict Lemma 5.4.

Next recall the sequences $\beta$ and $\gamma$ defined in (51). Thanks to (39),

$$\beta_{n+1} = \left( 1 - \frac{\theta}{n+1} \right) \beta_n + \theta \frac{(n+2)(n-3)}{n(n-1)(n+1)} \beta_{n-1},$$

which proves that $\beta_n$ has constant sign for large $n$, otherwise it would converge to 0 [and then $\Psi_\infty = 0$, which would contradict (i)]. Therefore, we can assume that $\beta_n \geq 0$ for large $n$ without loss of generality. Since

$$\beta_{n+1} - \beta_n = -\frac{\theta}{n+1}(\beta_n - \beta_{n-1}) - 6\theta \frac{\beta_{n-1}}{n(n-1)(n+1)},$$

then for any sufficiently large $n$, $0 \leq \beta_{n+1} < \max(\beta_n, \beta_{n-1})$, so that $(\beta_n)_n$ is bounded, and so is $(\gamma_n)_n$. Use (52) to show that

$$\beta_{n+1} - \beta_3 = -\theta \frac{n+2}{(n-1)(n-2)} \gamma_n - \theta \sum_{k=3}^{n-1} \frac{k+6}{k(k-1)(k-2)} \gamma_k,$$

and conclude that $(\beta_n)_n$ is convergent. $\square$

PROOF OF PROPOSITION 5.13. Recall $Z$ defined in Claim 1 and set $\varphi_2 := -(2\theta+3)/3c(\theta+3)$, as well as $V \in \mathcal{L}_3$,

$$V := \psi_{-2} e^{(-2)} + \psi_{-1} e^{(-1)} + \psi_1 e^{(1)} + \varphi_2 e^{(2)},$$

so that

(53) $$u^\varepsilon = -(c\Psi_\infty)^{-1}(Z - \Sigma e^{(2)}) + V + \frac{1}{\theta(\theta+3)} D^{(3)}.$$

By an elementary computation relying on (30) and (31), get

$$L'V = -\frac{1}{\theta(\theta+3)} \delta^{(3)} + e^{(-1)} - e^{(0)} - e^{(1)} + e^{(2)},$$

and conclude, thanks to Claim 1, that $L'u^\varepsilon = e^{(-1)} - e^{(0)}$.

To get the equivalent of $g_n^\varepsilon$ as $n \to \infty$, first recall (53) and observe that

$$\psi_{-2} + \psi_{-1} + \psi_1 + \varphi_2 = 0,$$

so that $(n^2 V_n)_n$ converges. Next consider $nZ_n$:

$$nZ_n = \sum_{k \geq 3} \frac{n\Psi_k}{n+k} = \sum_{k \geq 3} k\Psi_k \left( \frac{1}{k} - \frac{1}{n+k} \right)$$

$$= \Sigma - \sum_{k \geq 3} \frac{k\Psi_k}{n+k}$$

$$= \Sigma - \Psi_\infty \frac{\ln(n)}{n} + o\left(\frac{\ln(n)}{n}\right),$$

by a similar method as in the proof of Proposition 5.11. As a consequence,

$$nu_n^\varepsilon = -(c\Psi_\infty)^{-1}\left(nZ_n - \Sigma \frac{n}{n+2}\right) + nV_n + \frac{1}{\theta(\theta+3)} nD_n^{(3)}$$

$$= \frac{\ln(n)}{cn} + o\left(\frac{\ln(n)}{n}\right),$$

which ends the proof, since $g_n^\varepsilon = n^2 u_n^\varepsilon$. $\square$

## 6. Proof of the convergence to the TSS and the canonical diffusion.

6.1. *Preliminary result.* We start this section by stating and proving a technical proposition about the two type particular case without mutation ($\mu \equiv 0$) of the *GL*-population of Section 2.2. In the case where $\nu_0 = X_0 \delta_x + Y_0 \delta_y$, for any $t > 0$, $\nu_t = X_t \delta_x + Y_t \delta_y$. Then, the Markov process $(X_t, Y_t; t \geq 0)$ satisfies the following properties.

PROPOSITION 6.1.

(a) *Let $(\hat{X}_n, \hat{Y}_n; n \in \mathbb{N})$ denote the discrete-time Markov chain associated with $(X_t, Y_t; t \geq 0)$, and $\hat{T}$ denote the first hitting time of $\Omega_1 \cup \Omega_2$ by $(\hat{X}, \hat{Y})$. There is some positive constant $C$ independent of $x$ and $y$ such that*

$$(54) \qquad \mathbb{E}_{n,m}(\hat{T}) < C(n+m) \quad and \quad \mathbb{E}_{n,m}(\hat{T}^2) < C(n+m)^2.$$

(b) $\mathbb{E}_{n,m}(X_T^2 + Y_T^2) < C(n+m)^2$ *for some constant $C$ independent of $x$ and $y$.*

(c) *With the notation of Section 4.1, when the birth rates $b(x, y, n, m)$ are multiplied by a positive constant $a$, the fixation probabilities $u_{n,m}$ are continuous as a function of $a$.*

PROOF. (a) The process $(\hat{X}_n + \hat{Y}_n; n \in \mathbb{N})$ is dominated by the Markov chain $(\hat{Z}_n; n \in \mathbb{N})$ in $\mathbb{N}^*$ with initial state $k = X_0 + Y_0$ and transition probabilities

$$p_{ij} = \begin{cases} \bar{b}/[\bar{b} + \underline{c}(i-1)^\alpha], & \text{if } i \geq 1 \text{ and } j = i+1, \\ \underline{c}(i-1)^\alpha/[\bar{b} + \underline{c}(i-1)^\alpha], & \text{if } i \geq 2 \text{ and } j = i-1, \\ 0, & \text{otherwise.} \end{cases}$$

Let us denote by $P_k$ its law. Therefore, $\hat{T}$ is dominated by $\hat{S} := \inf\{n \geq 0, \hat{Z}_n = 1\}$ and it suffices to prove that $E_k(\hat{S}) \leq Ck$ and $E_k(\hat{S}^2) \leq Ck^2$ for some $C > 0$.

Let $(\tilde{U}_n; n \geq 0)$ be the discrete-time random walk on $\mathbb{Z}$ with right transition probability $1/3$ and left transition probability $2/3$. The law of $\tilde{U}$ conditional on $\tilde{U}_0 = k$ is denoted by $\tilde{P}_k$. Let $\tau$ be the first hitting time of 0 by $\tilde{U}$. For any $k \geq 0$, one can compute explicitly (see, e.g., [30]) that

$$\tilde{E}_k(\tau) = 3k, \qquad \tilde{E}_k(\tau^2) = 3k(3k+8) \quad \text{and} \quad \tilde{E}_k(\exp(\rho\tau)) = \exp(\sigma_\rho k)$$

for $0 \leq \rho \leq \ln(9/8)/2$ with $\exp(\sigma_\rho) = (1 - \sqrt{1 - 8\exp(2\rho)/9})3\exp(-\rho)/2$.

Now, let $k_0$ be large enough to have $\underline{c}(k_0 - 1)^\alpha > 2\bar{b}$. First, observe that any excursion of $\hat{Z}$ above $k_0 + 1$ is stochastically dominated by an excursion of the random walk $\tilde{U}$ above $k_0 + 1$. Second, let

$$\hat{S}' := \inf\{j \geq k_0 : \hat{Z}_j = 1, \hat{Z}_{j-1} = 2, \ldots, \hat{Z}_{j-k_0} = k_0 + 1\}.$$

Obviously, $\hat{S} \leq \hat{S}'$. Moreover, for any $n$ such that $\hat{Z}_n = k_0 + 1$, $\hat{S}' = n + k_0$ with probability $\beta := p_{k_0+1,k_0} \cdots p_{2,1} > 0$, and otherwise, $\hat{Z}_{n+k_0} \leq 2k_0 + 1$. Therefore, under $P_k$, $\hat{S}$ is dominated by

$$\tau_0 + \sum_{i=1}^{G}(k_0 + \tau_i),$$

where $G, \tau_0, \tau_1, \ldots$ are independent, $\tau_0$ has the law of the first hitting time of $k_0 + 1$ under $\tilde{P}_{k \vee (k_0+1)}$, $G$ is a geometric r.v. with parameter $\beta$, and the $\tau_i$'s are i.i.d. r.v. distributed as the first hitting time of $k_0 + 1$ under $\tilde{P}_{2k_0+1}$.

Therefore, if $\rho$ is small enough to have $(1-\beta)\exp(k_0(\rho+\sigma_\rho)) < 1$ (observe that $\sigma_\rho \to 0$ when $\rho \to 0$), then

$$E_k(e^{\rho \hat{S}}) \leq \tilde{E}_{(k-k_0-1)\vee 0}(e^{\rho \tau})\left(\sum_{j \geq 1} \beta(1-\beta)^{j-1}(e^{\rho k_0}\tilde{E}_{k_0}(e^{\rho \tau}))^j\right)$$

$$\leq \frac{\beta}{1-\beta}e^{\mu_\rho(k-k_0)}\sum_{j \geq 1}[(1-\beta)e^{k_0(\rho+\sigma_\rho)}]^j.$$

Therefore, there exists $C, \rho, \sigma > 0$ such that

(55) $$E_k(e^{\rho \hat{S}}) \leq Ce^{\sigma k}.$$

Moreover, for any $k \geq k_0$,

$$E_k(\hat{S}) \leq E_{k_0}(\hat{S}) + \tilde{E}_{k-k_0}(\tau) = E_{k_0}(\hat{S}) + 3(k - k_0)$$

and

$$E_k(\hat{S}^2) \leq 2E_{k_0}(\hat{S}^2) + 2\tilde{E}_{k-k_0}(\tau^2) = 2E_{k_0}(\hat{S}^2) + 6(k-k_0)(3(k-k_0)+8).$$

Since $E_{k_0}(\hat{S}^2) < \infty$ by (55), this ends the proof of (54).

(b) With the same notation as above, since $\hat{T}$ is the number of jumps of the process $(X, Y)$ that occurred on the time interval $[0, T]$, $X_T \leq X_0 + \hat{T}$ and $Y_T \leq Y_0 + \hat{T}$. Hence, $\mathbb{E}_{n,m}(X_T^2 + Y_T^2) \leq 2\mathbb{E}_{n,m}((n+m+\hat{T})^2)$ and the required bound follows from (a).

(c) Let us only denote the dependence of the fixation probability in $a$ by $u_{n,m}(a)$, and let us denote by $\pi_{ij}(a)$ the transitions (10) of the Markov chain $(\hat{X}, \hat{Y})$ when $b(x, y, n, m)$ and $b(y, x, n, m)$ are replaced by $ab(x, y, n, m)$ and $ab(y, x, n, m)$. It suffices to prove that $u_{n,m}(a)$ is continuous at $a = 1$.

We will use the notation $\pi_{i_1,\ldots,i_k}(a)$ for the product $\pi_{i_1 i_2}(a)\pi_{i_2 i_3}(a) \cdots \pi_{i_{k-1} i_k}(a)$ and $S_{(n,m)\to \Gamma}$ for the set of paths linking $(n, m)$ to a subset $\Gamma$ of $\mathbb{N}^2$ without hitting $\Omega_1 \cup \Omega_2$ before $\Gamma$, that is, the set of all $k$-tuples $(i_1, i_2, \ldots, i_k)$ for all $k \geq 1$ such that $i_1 = (n, m)$, $i_2, \ldots, i_{k-1} \in \mathbb{N}^2 \setminus (\Omega_1 \cup \Omega_2)$ and $i_k \in \Gamma$. Now,

$$u_{n,m}(a) = \sum_{k \geq 2} \sum_{(i_1,\ldots,i_k) \in S_{(n,m)\to \Omega_1}} \pi_{i_1,\ldots,i_k}(a),$$

so it is sufficient to prove that the previous series (in $k$) is uniformly convergent for $a$ in some neighborhood of 1.

Fix $i = (n, m)$ and $j \in \mathbb{N}^2$ such that $\|i - j\| = 1$, where $\|\cdot\|$ denotes the $L^1$ norm in $\mathbb{Z}^2$. It is elementary to check that $\pi_{ij}(a)/\pi_{ij}(1) \leq (a \vee 1)/(1 + (a-1)c_{n,m})$, where $c_{n,m} = B_{n,m}/(B_{n,m} + D_{n,m})$ with $B_{n,m} = nb(x, y, n, m) + mb(y, x, m, n)$ and $D_{n,m} = nd(x, y, n, m) + md(y, x, m, n)$. Because of the bounds we assumed on $b$ and $d$ (Definition 2.1), the $c_{n,m}$ are bounded by some constant $C$. Then, if $|a - 1|C < 1/2$, we finally get $\pi_{ij}(a) \leq 2(a \vee 1)\pi_{ij}(1)$, which implies the required uniform convergence. $\square$

6.2. *Proofs of Theorems* 3.1 *and* 3.2. Recall from Section 3 the definition of the stopping times $\tau_n$ ($n$th mutation time) and $\rho_n$ (first time after $\tau_n$ when the population becomes monomorphic) and of the random variables $V_n$ (the surviving trait at time $\rho_n$). Recall also the notation $b(x, n) \equiv b(x, n\delta_x)$ for the birth rate of an individual of trait $x$ in a monomorphic population made of $n$ individuals of trait $x$.

The proof relies on the following three lemmas. The first one states that there is no accumulation of mutations on the timescale $t/\gamma$. The second one gives the limiting laws of $\gamma\tau_1$ and of the population size at time $\tau_1$. The last one gives the behavior of $\rho_0$ and $V_0$ when the initial population is dimorphic.

LEMMA 6.2. *Fix $C, \eta > 0$. There is $\varepsilon > 0$ such that, for any $\gamma \in (0, 1)$,*

$$(56) \quad \mathbb{E}(\langle \nu_0^\gamma, \mathbf{1}\rangle) \leq C \quad \Longrightarrow \quad \forall t \geq 0 \quad \mathbb{P}\left(\exists n \in \mathbb{N}^* : \frac{t}{\gamma} \leq \tau_n \leq \frac{t+\varepsilon}{\gamma}\right) < \eta.$$

*Moreover, for any $\eta > 0$ and $t \geq 0$, there exists $n \in \mathbb{N}^*$ such that, for any $\gamma \in (0, 1)$,*

$$(57) \quad \mathbb{E}(\langle \nu_0^\gamma, \mathbf{1}\rangle) \leq C \quad \Longrightarrow \quad \mathbb{P}(\tau_n \leq t/\gamma) < \eta.$$

LEMMA 6.3. *Assume $\nu_0^\gamma = n\delta_x$ where $\mu(x) > 0$.*

(a) *As $\gamma \to 0$, the pair $(\gamma\tau_1, \langle \nu_{\tau_1-}^\gamma, \mathbf{1}\rangle)$ converges in law to a couple of independent random variables $(T, N)$, where $T$ is an exponentially distributed random variable with parameter $\beta(x)$ defined in (8), and the law of $N$ is obtained by $b(x, \cdot)$-size-biasing $\xi(x)$:*

$$(58) \quad \mathbb{P}(N = k) = \frac{kb(x, k)\mathbb{P}(\xi(x) = k)}{\mathbb{E}(\xi(x)b(x, \xi(x)))}.$$

(b) *For any $p \geq 1$, $\sup_{\gamma \in (0,1)} \mathbb{E}_{n\delta_x}^\gamma(\langle \nu_{\tau_1}, \mathbf{1}\rangle^p) < \infty$.*

LEMMA 6.4. *Assume $\nu_0^\gamma = n\delta_x + \delta_y$ (with $y \neq x$). Then:*

(a) *$\gamma\rho_0 \to 0$ in probability and $\mathbb{P}^\gamma(\rho_0 < \tau_1) \to 1$ as $\gamma \to 0$.*

(b) $\sup_{\gamma \in (0,1)} \mathbb{E}^\gamma(\langle \nu_{\rho_0}^\gamma, \mathbf{1} \rangle^2 \mathbf{1}_{\{\rho_0 < \tau_1\}}) < \infty$.
(c) $\lim_{\gamma \to 0} \mathbb{P}^\gamma(V_0 = y) = 1 - \lim_{\gamma \to 0} \mathbb{P}^\gamma(V_0 = x) = u_{n,1}(x,y)$.

PROOF OF LEMMA 6.2. Fix $C > 0$ and assume $\mathbb{E}(\langle \nu_0^\gamma, \mathbf{1} \rangle) \leq C$. By Proposition 2.3, there exists a constant $C'$ such that $\mathbb{E}(\langle \nu_t^\gamma, \mathbf{1} \rangle) \leq C'$ for any $t \geq 0$ and $\gamma > 0$. Therefore, it is sufficient to show (56) for $t = 0$.

Now, when the total population size is $n$, the total mutation rate in the population is bounded by $\gamma \bar{b} n$, so that the number of mutations $M_t$ between times 0 and $t$ is dominated by a Poisson point process with intensity $\gamma \bar{b} \langle \nu_s, \mathbf{1} \rangle \, ds$. Therefore,

$$\mathbb{P}(M_{\varepsilon/\gamma} \geq 1) \leq \mathbb{E}(M_{\varepsilon/\gamma}) \leq \gamma \bar{b} \int_0^{\varepsilon/\gamma} \mathbb{E}(\langle \nu_s^\gamma, \mathbf{1} \rangle) \, ds \leq \varepsilon \bar{b} C',$$

which concludes the proof of (56).

Similarly, for $t \geq 0$, $\mathbb{P}(M_{t/\gamma} \geq n) \leq t \bar{b} C'/n$, which implies (57). □

PROOF OF LEMMA 6.3. Fix $\gamma \in (0,1)$ and assume that $\nu_0^\gamma = n \delta_x$. Recall that $\tau_1$ is the first mutation time. Notice that $(\nu_t^\gamma; t < \tau_1)$ is distributed as $(X_t^\gamma \delta_x; t < \tau)$, where $X^\gamma$ is a birth-and-death process with initial state $X_0^\gamma = n$ and transition rates $(1 - \gamma \mu(x)) i b(x, i)$ from $i$ to $i+1$ and $i d(x, i)$ from $i$ to $i-1$, and $\tau$ is the first point of a Poisson point process with inhomogeneous intensity $g^\gamma(X_t^\gamma) := \gamma \mu(x) X_t^\gamma b(x, X_t^\gamma)$ (depending on $X^\gamma$ solely through its intensity). Therefore, for any bounded function $f : \mathbb{N}^* \to \mathbb{R}$ and for any $t \geq 0$,

$$(59) \quad \mathbb{E}^\gamma(f(\langle \nu_{t/\gamma}^\gamma, \mathbf{1} \rangle); \gamma \tau_1 > t) = \mathbb{E}\left[ f(X_{t/\gamma}^\gamma) \exp\left(-\int_0^{t/\gamma} g^\gamma(X_s^\gamma) \, ds\right) \right]$$

and

$$\mathbb{E}^\gamma(f(\langle \nu_{\tau_1-}^\gamma, \mathbf{1} \rangle); \gamma \tau_1 \leq t) = \mathbb{E}\left( \int_0^{t/\gamma} f(X_s^\gamma) g(X_s^\gamma) e^{-\int_0^s g(X_u^\gamma) du} \, ds \right),$$

which yields, after a change of variable,

$$(60) \quad \begin{aligned} &\mathbb{E}^\gamma(f(\langle \nu_{\tau_1-}^\gamma, \mathbf{1} \rangle) \mathbf{1}_{\{\gamma \tau_1 \leq t\}}) \\ &= \mu(x) \int_0^t \mathbb{E}(f(X_{s/\gamma}^\gamma) X_{s/\gamma}^\gamma b(x, X_{s/\gamma}^\gamma) e^{-\gamma \mu(x) \int_0^{s/\gamma} X_u^\gamma b(x, X_u^\gamma) du}) \, ds. \end{aligned}$$

Now, since the individual birth rates of $X^\gamma$ decrease with $\gamma$, all the processes $X^\gamma$ can be constructed on a same space in such a way that, for $0 \leq \gamma \leq \gamma' \leq 1$ and $t \geq 0$, $X_t^{\gamma'} \leq X_t^\gamma \leq X_t^0$.

To compute the limit of (60) when $\gamma \to 0$, let us first prove that

$$(61) \quad \lim_{\gamma \to 0} \left| \int_0^t \mathbb{E}(g_1(X_{s/\gamma}^\gamma) e^{-\gamma \int_0^{s/\gamma} g_2(X_u^\gamma) du}) \, ds - \int_0^t \mathbb{E}(g_1(X_{s/\gamma}^0) e^{-\gamma \int_0^{s/\gamma} g_2(X_u^0) du}) \, ds \right| = 0$$

for any functions $g_1$ and $g_2$ on $\mathbb{N}$ such that $|g_1(x)| \leq Cx$ and $0 \leq g_2(x) \leq Cx$. For $M > 0$, this quantity is bounded by

$$\int_0^t \mathbb{E}|g_1(X^\gamma_{s/\gamma}) - g_1(X^0_{s/\gamma})|\, ds + C \int_0^t \gamma \int_0^{s/\gamma} \mathbb{E}(X^\gamma_{s/\gamma}|g_2(X^\gamma_u) - g_2(X^0_u)|)\, du\, ds$$

$$\leq 2C\gamma \int_0^{t/\gamma} \mathbb{E}(X^0_u; X^0_u \neq X^\gamma_u)\, du$$

$$+ 2C^2 \int_0^t \gamma \int_0^{s/\gamma} \mathbb{E}(X^\gamma_{s/\gamma} X^0_u; X^0_u \neq X^\gamma_u)\, du\, ds$$

$$\leq +2C\gamma \int_0^{t/\gamma} [\mathbb{E}(X^0_u; X^0_u > M) + M\mathbb{P}(X^0_u \leq M, X^0_u \neq X^\gamma_u)]\, du$$

$$+ 2C^2 \int_0^t \gamma \int_0^{s/\gamma} [\mathbb{E}(X^\gamma_{s/\gamma} X^0_u; X^0_u > M)$$

$$+ M\mathbb{E}(X^\gamma_{s/\gamma}; X^0_u \leq M, X^0_u \neq X^\gamma_u)]\, du\, ds.$$

Therefore, by Proposition 2.3, using the Cauchy–Schwarz inequality to bound the term involving $\mathbb{E}(X^\gamma_{s/\gamma} X^0_u; X^0_u > M)$ and distinguishing between $X^\gamma_{s/\gamma} \leq M$ and $X^\gamma_{s/\gamma} > M$ in the term involving $\mathbb{E}(X^\gamma_{s/\gamma}; X^0_u \leq M, X^0_u \neq X^\gamma_u)$, it is sufficient to prove that

$$\lim_{\gamma \to 0} \gamma \int_0^{t/\gamma} \mathbb{P}(X^0_u \leq M, X^0_u \neq X^\gamma_u)\, du = 0$$

or, equivalently, that $\gamma$ times the expected time length between 0 and $t/\gamma$ where $X^0_u \leq M$ and $X^0_u \neq X^\gamma_u$ goes to 0.

Since the difference between the birth rates of $X^0_u$ and $X^\gamma_u$ when $X^0_u = X^\gamma_u$ is less than $\gamma \bar{b} X^0_u$, any time when the two processes can start to differ belongs to the set of times when a point of a Poisson point process on $\mathbb{R}^2_+$ with intensity $\gamma \bar{b}\, du\, ds$ independent of $X^0$ lies below the curve $(t, X^0_t)_{t \geq 0}$. For each of these points, the time length where the two processes differ is dominated by the first hitting time of 1 by $X^0$ (at this time, the two processes are necessarily equal). Since, moreover, we only have to consider the time intervals where $X^0_u \leq M$, all these time lengths are dominated by independent realizations of the hitting time of 1 by $X^0$ starting from $M$.

Let $n_\gamma$ denote the number of points of the previous Poisson point process below $X^0$ before time $t/\gamma$ and let $(R_k)$ be a sequence of r.v. independent of $n_\gamma$ distributed as the hitting time of 1 by $X^0$ starting from $M$. Then

$$\gamma \int_0^{t/\gamma} \mathbb{P}(X^0_u \leq M, X^0_u \neq X^\gamma_u)\, du$$

$$\leq \gamma \mathbb{E}\left(\sum_{k=1}^{n_\gamma} R_k\right) = \gamma \mathbb{E}(R_1)\mathbb{E}(n_\gamma) = \gamma \mathbb{E}(R_1)\bar{b}\gamma \int_0^{t/\gamma} \mathbb{E}(X^0_u)\, du,$$

which goes to 0 as $\gamma \to 0$ by Proposition 2.3.

Therefore, (61) is proved. Combined with (60), we get

$$\lim_{\gamma \to 0} \mathbb{E}^\gamma(f(\langle \nu^\gamma_{\tau_1-}, \mathbf{1}\rangle)\mathbf{1}_{\{\gamma\tau_1 \leq t\}})$$
$$= \lim_{\gamma \to 0} \mu(x) \int_0^t \mathbb{E}(f(X^0_{s/\gamma})X^0_{s/\gamma}b(x, X^0_{s/\gamma})e^{-\gamma\mu(x)\int_0^{s/\gamma} X^0_u b(x,X^0_u)du})\,ds.$$

Since $X^0$ is exactly the positive-recurrent Markov chain mentioned at the end of Proposition 2.2, apply the ergodic theorem to get that $\gamma \int_0^{s/\gamma} X^0_u b(x, X^0_u)\,du \to s\mathbb{E}(\xi(x)b(x,\xi(x)))$ a.s. as $\gamma \to 0$. Since $\sup_{t \geq 0} \mathbb{E}((X^0_t)^2) < \infty$ and by Lebesgue's theorem, we finally get

$$\lim_{\gamma \to 0} \mathbb{E}^\gamma(f(\langle \nu^\gamma_{\tau_1-}, \mathbf{1}\rangle)\mathbf{1}_{\{\gamma\tau_1 \leq t\}})$$
$$= \mu(x) \int_0^t \mathbb{E}(f(\xi(x))\xi(x)b(x,\xi(x))e^{-\mu(x)\mathbb{E}[\xi(x)b(x,\xi(x))]s})$$
$$= \frac{\mathbb{E}(f(\xi(x))\xi(x)b(x,\xi(x)))}{\mathbb{E}(\xi(x)b(x,\xi(x)))} \int_0^t \beta(x)e^{-\beta(x)s}\,ds,$$

which completes the proof of Lemma 6.3(a).

Lemma 6.3(b) can be obtained by taking $f(x) = x^p \wedge K$ in (60), then letting first $K$ go to infinity and next $t$ to infinity. Then, we get that

$$\mathbb{E}^\gamma_{n\delta_x}(\langle \nu_{\tau_1}, \mathbf{1}\rangle^p) \leq \bar{b} \int_0^\infty \mathbb{E}\left[(X^\gamma_{s/\gamma})^{p+1}\exp\left(-\gamma\mu(x)\int_0^{s/\gamma} X^\gamma_u b(x, X^\gamma_u)\,du\right)\right]ds$$
$$\leq \bar{b} \int_0^\infty \mathbb{E}[(X^\gamma_{s/\gamma})^{p+1}\exp(-\gamma\mu(x)b(x,1)L_{s/\gamma})]\,ds,$$

where $L_t = \int_0^t \mathbf{1}_{\{X^0_u = 1\}}\,du$. By the ergodic theorem for $X^0$, $L_t/t$ converges a.s. as $t \to \infty$ to a positive nonrandom limit $l$. However, we need a finer result to conclude. Fix $\lambda > 0$. Distinguishing between $L_{s/\gamma} \leq \lambda s/\gamma$ and $L_{s/\gamma} > \lambda s/\gamma$ and using the Cauchy–Schwarz inequality and Proposition 2.3, one can bound from above the last displayed integral by a constant times

$$\int_0^\infty [\mathbb{P}(L_{s/\gamma} \leq \lambda s/\gamma)^{1/2} + \exp(-\mu(x)b(x,1)\lambda s)]\,ds.$$

Therefore, it suffices to prove that there exist $\lambda, \lambda', C > 0$ such that $\mathbb{P}(L_t \leq \lambda t) \leq Ce^{-\lambda' t}$ for any $t \geq 0$.

Now, define recursively $t_0 = 0$, and for $i \geq 1$, $s_i = \inf\{s \geq t_{i-1} : X^0_s = 1\}$, $t_i = \inf\{t \geq s_i : X^0_s = 2\}$. Then for $i \geq 1$, set $T_i := t_i - s_i$ and $S_i := s_i - t_{i-1}$. By the strong Markov property, all these r.v. are independent, and more specifically, $(T_i)_{i \geq 1}$ are i.i.d. exponential r.v. with parameter $b(x, 1)$ (the jump rate of $X^0$ from state 1), $(S_i)_{i \geq 2}$ are i.i.d. r.v. distributed as the hitting time of 1 by the process $X^0$ started at 2, whereas $S_1$ is the hitting time of 1 by the process $X^0$ started at $n$ (remember that $n$ has been defined by $\nu_0 = n\delta_x$). Then, for any $\rho, \sigma > 0$, using Chebyshev's exponential inequality to get the last line,

$$\mathbb{P}(L_t \leq \lambda t) \leq \mathbb{P}\left(\exists k \geq 1 : \sum_{i=1}^{k} T_i \leq \lambda t \text{ and } \sum_{i=1}^{k+1}(S_i + T_i) \geq t\right)$$

$$\leq \mathbb{P}\left(\exists k \geq 1 : \sum_{i=1}^{k} T_i \leq \lambda t \text{ and } T_{k+1} + \sum_{i=1}^{k+1} S_i \geq (1-\lambda)t\right)$$

$$\leq \sum_{k=1}^{\infty} \mathbb{P}\left(\sum_{i=1}^{k} T_i \leq \lambda t\right) \mathbb{P}\left(T_{k+1} + \sum_{i=1}^{k+1} S_i \geq (1-\lambda)t\right)$$

$$\leq \mathbb{E}(e^{\sigma(T_1+S_1)}) e^{[\rho\lambda - \sigma(1-\lambda)]t} \sum_{k=1}^{\infty} [\mathbb{E}(e^{-\rho T_1})\mathbb{E}(e^{\sigma S_2})]^k.$$

Observe that $\mathbb{E}(\exp(\sigma T_1)) = b(x,1)/(b(x,1) - \sigma)$ if $\sigma < b(x,1)$. Therefore, if we can prove that there exists $\sigma > 0$ such that $\mathbb{E}(\exp(\sigma S_1)) < \infty$ [and thus $\mathbb{E}(\exp(\sigma S_2)) < \infty$], then $\rho$ can be chosen large enough to have $\mathbb{E}(\exp(-\rho T_1)) \times \mathbb{E}(\exp(\sigma S_2)) < 1$, and next $\lambda > 0$ can be chosen small enough to have $\rho\lambda - \sigma(1-\lambda) < 0$. This would end the proof of Lemma 6.3(b).

Therefore, it only remains to prove that there exists $\sigma > 0$ such that $\mathbb{E}(\exp(\sigma \times S_1)) < \infty$. This can be done as follows. Let $\hat{X}^0$ be the discrete-time Markov chain associated with $X^0$ and let $U_k$ be the holding time of $X^0$ in the state $\hat{X}^0_k$. Let also $\hat{S}_1$ be the first integer $k$ such that $\hat{X}^0_k = 1$. Then, for any $k$, $U_k$ is dominated by an exponential r.v. with parameter $\kappa := b(x,1) \vee \underline{c}$, which is a lower bound for the jump rates of $X^0$. Therefore, if $\sigma < \kappa$ and $(R_i)_{i \geq 1}$ denote i.i.d. exponential r.v. with parameter $\kappa$,

$$\mathbb{E}(e^{\sigma S_1}) \leq \sum_{k=0}^{\infty} \mathbb{E}\left(\exp\left(\sigma \sum_{i=1}^{k} R_i\right)\right) \mathbb{P}(\hat{S}_1 = k)$$

$$= \mathbb{E}\left[\left(\frac{\kappa}{\kappa - \sigma}\right)^{\hat{S}_1}\right],$$

which has already been proved to be finite for small enough $\sigma$ in (55). $\square$

PROOF OF LEMMA 6.4. Before the first mutation, $\nu_t^\gamma = X_t^\gamma \delta_x + Y_t^\gamma \delta_y$, where $(X_t^\gamma, Y_t^\gamma)$ is a two-type $GL$-population with transition rates

$$\begin{cases} (1 - \gamma\mu(x))nb(x,y,n,m), & \text{from } (n,m) \text{ to } (n+1,m), \\ (1 - \gamma\mu(y))mb(y,x,m,n), & \text{from } (n,m) \text{ to } (n,m+1), \\ nd(x,y,n,m), & \text{from } (n,m) \text{ to } (n-1,m), \\ mb(y,x,m,n), & \text{from } (n,m) \text{ to } (n,m-1), \\ 0, & \text{otherwise.} \end{cases}$$

On the event $\{\tau_1 > \rho_0\}$, $V_0 = y$ if and only if there is fixation in this process, $V_0 = x$ otherwise, and $\rho_0$ equals the fixation time $T$ (see Section 2.2.2).

Now, by Lemma 6.2, for any $\eta > 0$, there exists $\varepsilon > 0$ such that $\mathbb{P}(\tau_1 > \varepsilon/\gamma) \geq 1 - \eta$. Since $\mathbb{P}_{n,1}(T < \infty) = 1$, this implies easily (a). It is then elementary to deduce from Proposition 6.1(c) that (c) holds. Finally, (b) follows from the observation that

$$\mathbb{E}_{n\delta_x + \delta_y}(\langle \nu^\gamma_{\rho_0}, \mathbf{1} \rangle^2 \mathbf{1}_{\{\rho_0 < \tau_1\}}) \leq \mathbb{E}_{n,1}((X^\gamma_T)^2 + (Y^\gamma_T)^2)$$

and from Proposition 6.1(b). □

PROOF OF THEOREM 3.1. Observe that the generator $A$ of the process $S$, defined in (7), can be written as

(62) $$A\varphi(x) = \int_{\mathcal{X}} (\varphi(x+h) - \varphi(x)) \beta(x) \kappa(x, dh),$$

where $\beta(x)$ has been defined in (8) and where $\kappa(x, dh)$ is the *probability* measure on $\mathcal{X} - x$ defined by

(63) $$\kappa(x, dh) = \sum_{n=1}^{\infty} u_{n,1}(x, x+h) \frac{nb(x,n)\mathbb{P}(\xi(x) = n)}{\mathbb{E}(\xi(x)b(x,\xi(x)))} M(x, dh)$$
$$+ \left( \int_{\mathbb{R}^k} \sum_{n=1}^{\infty} (1 - u_{n,1}(x, x+y)) \right.$$
$$\left. \times \frac{nb(x,n)\mathbb{P}(\xi(x) = n)}{\mathbb{E}(\xi(x)b(x,\xi(x)))} M(x, dy) \right) \delta_0(dh).$$

This means that the TSS model $S$ with initial state $x$ can be constructed as follows: let $(U(k), k = 0, 1, 2, \ldots)$ be a Markov chain in $\mathcal{X}$ with initial state $x$ and with transition kernel $\kappa(x, dh)$, and let $(P(t), t \geq 0)$ be an independent standard Poisson process. Then, the process $(S_t, t \geq 0)$ defined by

$$S_t = U \circ P\left( \int_0^t \beta(S_s) \, ds \right)$$

is a Markov process with infinitesimal generator (62), cf. [10], Chapter 6. Let $(J_n)_{n \geq 1}$ denote the sequence of jump times of the Poisson process $P$ and define $(T_n)_{n \geq 1}$ by $J_n = \int_0^{T_n} \beta(S_s) \, ds$ or $T_n = \infty$ if $\int_0^\infty \beta(S_s) \, ds < J_n$. Observe that any jump of the process $S$ occurs at some time $T_n$, but that all $T_n$ may not be effective jump times for $S$, because of the Dirac mass at 0 appearing in (63). As will appear below, the sequence $(T_n)$ can be interpreted as the sequence of mutation times in the limit process. Whether an effective jump occurs at time $T_n$ or not then corresponds to the fixation or extinction of the mutant.

Let $\mathbf{P}_x$ denote the law of $\zeta_t$, as defined in the theorem, conditional on $\mathrm{Supp}(\zeta_0) = S_0 = x$. Fix $t > 0$, $m \in \mathbb{N}^*$, $x \in \mathcal{X}$ and a measurable subset $\Gamma$ of $\mathcal{X}$. Under $\mathbf{P}_x$, $T_1$ and $S_{T_1}$ are independent, $T_1$ is an exponential random variable with parameter $\beta(x)$, and $S_{T_1}$ has law $\kappa(x, \cdot)$. Therefore, for any $n \geq 1$, applying the strong Markov property to the process $S$ at time $T_1$ in the third line,

$$
\begin{aligned}
&\mathbf{P}_x(T_n \leq t < T_{n+1}, \exists z \in \Gamma : \zeta_t = m\delta_z) \\
&\quad = \mathbf{E}_x(\mathbf{1}_{\{S_t \in \Gamma\}} \mathbf{1}_{\{T_n \leq t < T_{n+1}\}} \mathbb{P}(\xi(S_t) = m)) \\
&\quad = \int_0^t \beta(x) e^{-\beta(x)s} \int_{\mathbb{R}^k} \mathbf{P}_{x+h}(T_{n-1} \leq t - s < T_n, \exists z \in \Gamma : \zeta_{t-s} = m\delta_z) \\
&\qquad\qquad\qquad\qquad \times \kappa(x, dh)\, ds.
\end{aligned}
\tag{64}
$$

Moreover,

$$
\mathbf{P}_x(0 \leq t < T_1, \exists z \in \Gamma : \zeta_t = m\delta_z) = \mathbf{1}_{\{x \in \Gamma\}} e^{-\beta(x)t} \mathbb{P}(\xi(x) = m). \tag{65}
$$

These two relations characterize the one-dimensional laws of the process $\zeta$. The idea of our proof is to show that the same relations hold when we replace $T_n$ by $\tau_n$ and the support of $\zeta_t$ by the support of $\nu^\gamma_{t/\gamma}$ (when it is a singleton) in the limit $\gamma \to 0$.

More precisely, let us define for any $\nu_0 \in \mathcal{M}$ and $n \in \mathbb{N}$,

$$
p_n^\gamma(t, \Gamma, m, \nu_0) := \mathbb{P}^\gamma_{\nu_0}\left(\rho_n \leq \frac{t}{\gamma} < \tau_{n+1}, \exists z \in \Gamma : \nu_{t/\gamma} = m\delta_z\right).
$$

We will prove the following lemma after the end of this proof.

LEMMA 6.5. *For any $x \in \mathcal{X}$, $m, k \geq 1$, $n \geq 0$, $t > 0$ and any measurable subset $\Gamma$ of $\mathcal{X}$, $p_n(t, \Gamma, m, x) := \lim_{\gamma \to 0} p_n^\gamma(t, \Gamma, m, k\delta_x)$ exists, is independent of $k$ and satisfies*

$$
p_0(t, \Gamma, m, x) = \mathbf{1}_{\{x \in \Gamma\}} e^{-\beta(x)t} \mathbb{P}(\xi(x) = m) \tag{66}
$$

*and*

$$
\forall n \geq 1 \quad p_n(t, \Gamma, m, x) = \int_0^t \beta(x) e^{-\beta(x)s} \int_{\mathbb{R}^k} p_{n-1}(t - s, \Gamma, m, x + h) \\ \times \kappa(x, dh)\, ds. \tag{67}
$$

Comparing (64) and (65) with (66) and (67), this lemma implies that $p_n(t, \Gamma, m, x) = \mathbf{P}_x(T_n \leq t < T_{n+1}, \exists z \in \Gamma : \zeta_t = m\delta_z)$.

Recall that $\nu_0^\gamma = N_0^\gamma \delta_x$ with $\sup_{\gamma \in (0,1)} \mathbb{E}((N_0^\gamma)^p) < \infty$ for some $p > 1$. By Proposition 2.3, $\sup_{\gamma \in (0,1)} \sup_{t \geq 0} \mathbb{E}(\langle \nu_t^\gamma, \mathbf{1} \rangle^p) < +\infty$ and

$$|\mathbb{P}_{\nu_0^\gamma}^\gamma(\exists z \in \Gamma : \nu_{t/\gamma}^\gamma = m\delta_z) - \mathbf{P}_x(\exists z \in \Gamma : \zeta_t = m\delta_z)|$$

$$\leq 2\mathbb{P}(N_0^\gamma > M) + \sum_{k=1}^{M} \left| \sum_{n=0}^{\infty} (p_n^\gamma(t, \Gamma, m, k\delta_x) - p_n(t, \Gamma, m, x)) \right| \mathbb{P}(N_0^\gamma = k).$$

Because of Lemma 6.2 (57), the quantity inside the absolute value in the r.h.s. of this equation vanishes as $\gamma \to 0$. Thus,

(68) $$\lim_{\gamma \to 0} \mathbb{P}_{\nu_0^\gamma}^\gamma(\exists z \in \Gamma : \nu_{t/\gamma}^\gamma = m\delta_z) = \mathbf{P}_x(\exists z \in \Gamma : \zeta_t = m\delta_z).$$

Taking $\Gamma = \mathcal{X}$ and summing this relation over $m \in \mathbb{N}^*$, Fatou's lemma yields

(69) $$\lim_{\gamma \to 0} \mathbb{P}_{\nu_0^\gamma}^\gamma(\mathrm{Supp}(\nu_{t/\gamma}^\gamma) \text{ is a singleton}) = 1.$$

Now, consider a bounded measurable $f : \mathcal{M} \to \mathbb{R}$ such that $f(\nu) = 0$ if $\langle \nu, \mathbf{1} \rangle \neq m \in \mathbb{N}^*$ and define the function $\hat{f} : \mathcal{X} \to \mathbb{R}$ by $\hat{f}(x) = f(m\delta_x)$. Then, it follows from (68) and (69) that

$$\lim_{\gamma \to 0} \mathbb{E}_{\nu_0^\gamma}^\gamma(f(\nu_{t/\gamma}^\gamma))$$

(70) $$= \lim_{\gamma \to 0} \mathbb{E}_{\nu_0^\gamma}^\gamma(\hat{f}(\mathrm{Supp}(\nu_{t/\gamma}^\gamma)); \mathrm{Supp}(\nu_{t/\gamma}^\gamma) \text{ is a singleton and } \langle \nu_{t/\gamma}^\gamma, \mathbf{1} \rangle = m)$$

$$= \mathbf{E}_x(\hat{f}(S_t); N_t = m) = \mathbf{E}_x(f(\zeta_t)).$$

This equality generalizes to any bounded measurable $f : \mathcal{M} \to \mathbb{R}$ using once again that $\sup_{\gamma \in (0,1)} \mathbb{E}(\langle \nu_{t/\gamma}^\gamma, \mathbf{1} \rangle^p) < +\infty$. This completes the proof of Theorem 3.1 for one-dimensional distributions.

The extension to finite dimensional marginals can be proved exactly in the same fashion. □

PROOF OF LEMMA 6.5. Recall $\beta(x) = \mu(x)\mathbb{E}(\xi(x)b(x,\xi(x)))$. The result is trivial when the mutation rate $\mu(x)$ is 0. Let us assume that $\mu(x) \neq 0$. We will prove this lemma by induction over $n \geq 0$.

Fix $x \in \mathcal{X}$, $m, k \geq 1$ and $t > 0$. First, we have already proved in (59) that

$$p_0^\gamma(t, \Gamma, m, k\delta_x) = \mathbf{1}_{\{x \in \Gamma\}} \mathbb{P}_{k\delta_x}^\gamma(\langle \nu_{t/\gamma}^\gamma, \mathbf{1} \rangle = m, \gamma \tau_1 > t)$$

$$= \mathbf{1}_{\{x \in \Gamma\}} \mathbb{E}\left[ \mathbf{1}_{\{X_{t/\gamma}^\gamma = m\}} \exp\left( -\gamma \mu(x) \int_0^{t/\gamma} X_{s-}^\gamma b(x, X_{s-}^\gamma) \, ds \right) \right],$$

where $X^\gamma$ has been defined in the proof of Lemma 6.3. Using (61) and the ergodic theorem for positive-recurrent Markov chains, we get

$$\lim_{\gamma \to 0} p_0^\gamma(t, \Gamma, m, k\delta_x) = \mathbf{1}_{\{x \in \Gamma\}} \mathbb{P}(\xi(x) = m) \exp(-\mu(x)\mathbb{E}(\xi(x)b(x,\xi(x)))t),$$

which entails (66).

Then, fix $n \geq 1$ and assume that $\lim_{\gamma \to 0} p_{n-1}^\gamma(t, \Gamma, m, k\delta_x) = p_n(t, \Gamma, m, x)$ for any $t > 0, \Gamma \subset \mathcal{X}, x \in \mathcal{X}$ and $m, k \geq 1$. Applying the strong Markov property to the process $\nu^\gamma$ at time $\tau_1$, and using the fact that the mutant trait at this time is $x + U$, where $U$ has law $M(x, dh)$ and is independent of $\nu^\gamma_{\tau_1-}$, we get

$$
\begin{aligned}
(71) \quad & p_n^\gamma(t, \Gamma, m, k\delta_x) \\
& = \int_{\mathbb{R}^k} \mathbb{E}_{k\delta_x}^\gamma [\mathbf{1}_{\{\gamma\tau_1 \leq t\}} p_{n-1}^\gamma(t - \gamma\tau_1, \Gamma, m, \langle \nu^\gamma_{\tau_1-}, \mathbf{1}\rangle \delta_x + \delta_{x+h})] M(x, dh).
\end{aligned}
$$

Now, we want to apply the strong Markov property to $\nu^\gamma$ at time $\rho_0$ to compute the quantity $p_{n-1}^\gamma(s, \Gamma, m, l\delta_x + \delta_y)$ appearing inside the expectation in the last formula. For $K > 0$, distinguishing between the cases where $\rho_0 > \tau_1$, $\langle \nu_{\rho_0}, \mathbf{1}\rangle > K$, $V_0 = x$ and $V_0 = y$ yields

$$
\begin{aligned}
& p_{n-1}^\gamma(s, \Gamma, m, l\delta_x + \delta_y) \\
& = \mathbb{E}_{l\delta_x + \delta_y}^\gamma [\mathbf{1}_{\{\rho_0 < \tau_1, \langle \nu^\gamma_{\rho_0}, \mathbf{1}\rangle \leq K, V_0 = x\}} p_{n-1}^\gamma(s - \gamma\rho_0, \Gamma, m, \langle \nu^\gamma_{\rho_0}, \mathbf{1}\rangle \delta_x)] \\
& \quad + \mathbb{E}_{l\delta_x + \delta_y}^\gamma [\mathbf{1}_{\{\rho_0 < \tau_1, \langle \nu^\gamma_{\rho_0}, \mathbf{1}\rangle \leq K, V_0 = y\}} p_{n-1}^\gamma(s - \gamma\rho_0, \Gamma, m, \langle \nu^\gamma_{\rho_0}, \mathbf{1}\rangle \delta_y)] \\
& \quad + \mathbb{P}_{l\delta_x + \delta_y}^\gamma(\{\rho_0 \geq \tau_1\} \cap E) + \mathbb{P}_{l\delta_x + \delta_y}^\gamma(\{\rho_0 < \tau_1\} \cap \{\langle \nu^\gamma_{\rho_0}, \mathbf{1}\rangle > K\} \cap E),
\end{aligned}
$$

where

$$E = \{\rho_{n-1} \leq s/\gamma < \tau_n, \exists z \in \Gamma : \nu_{s/\gamma} = m\delta_z\}.$$

The third term of the r.h.s. vanishes as $\gamma \to 0$ because of Lemma 6.4(a) and the last term vanishes as $K \to +\infty$ uniformly for $\gamma \in (0,1)$ because of Lemma 6.4(b).

As for the first two terms, assume that $p_{n-1}^\gamma(t, \Gamma, m, k\delta_x)$ converges to $p_{n-1}(t, \Gamma, m, x)$ as in the statement of Lemma 6.5. As a consequence of Lemma 6.2 (56), for any $t > 0$, the function $s \mapsto p_{n-1}^\gamma(s, \Gamma, m, k\delta_x)$ is uniformly continuous on $[0, t]$, uniformly in $\gamma$. Combining this observation with Lemma 6.4(c),

$$
\begin{aligned}
(72) \quad & \lim_{\gamma \to 0} p_{n-1}^\gamma(s, \Gamma, m, l\delta_x + \delta_y) \\
& = u_{l,1}(x, y) p_{n-1}(s, \Gamma, m, y) + (1 - u_{l,1}(x, y)) p_{n-1}(s, \Gamma, m, x).
\end{aligned}
$$

This uniform continuity argument also applies to $s \mapsto p_{n-1}^\gamma(s, \Gamma, m, l\delta_x + \delta_y)$, so that the convergence in (72) is uniform in $s \in [0, t]$ and $l \in \{1, \ldots, L\}$, for fixed $L \geq 1$. Therefore, distinguishing between $\langle \nu^\gamma_{\tau_1-}, \mathbf{1}\rangle \leq L$ and $\langle \nu^\gamma_{\tau_1-}, \mathbf{1}\rangle > L$, we can combine Lemma 6.3(a) and (b) to get

$$\lim_{\gamma \to 0} \mathbb{E}_{k\delta_x}^\gamma [\mathbf{1}_{\{\gamma\tau_1 \leq t\}} p_{n-1}^\gamma(t - \gamma\tau_1, \Gamma, m, \langle \nu^\gamma_{\tau_1-}, \mathbf{1}\rangle \delta_x + \delta_y)]$$

$$= \int_0^t ds\, \beta(x) e^{-\beta(x)s}$$
$$\times \sum_{l=1}^\infty \frac{lb(x,l)\mathbb{P}(\xi(x)=l)}{\mathbb{E}(\xi(x)b(x,\xi(x)))}$$
$$\times [u_{l,1}(x,y)p_{n-1}(t-s,\Gamma,m,y)$$
$$+ (1-u_{l,1}(x,y))p_{n-1}(t-s,\Gamma,m,x)].$$

Finally, using Lebesgue's theorem, this limit applies inside the integral in (71), which gives exactly (67) and ends the proof of Lemma 6.5. □

PROOF OF THEOREM 3.2. Since the limiting law of the process $(S_t^\gamma; t \geq 0)$ is characterized by its finite-dimensional distributions, obtained in Theorem 3.1, we only have to show the tightness of their laws. Fix $T > 0$. By Ascoli's theorem for cdlg functions (see, e.g., [2]), we have to show that, for any $\varepsilon, \eta > 0$, there is $\delta > 0$ such that

(73) $$\limsup_{\gamma \to 0} \mathbb{P}(\omega'(S^\gamma, \delta) > \eta) < \varepsilon,$$

where the modulus of continuity $\omega'$ is defined as

$$\omega'(f, \delta) = \inf\left\{\max_{0 \leq i \leq r-1} \omega(f, [t_i, t_{i+1}))\right\},$$

where the infimum is taken over all finite partitions $0 = t_0 < t_1 < \cdots < t_r = T$ of $[0,T]$ such that $t_{i+1} - t_i > \delta$ for any $0 \leq i \leq r-1$, and where $\omega$ is the usual modulus of continuity

$$\omega(f, I) = \sup_{s,t \in I} \|f(t) - f(s)\|.$$

Now, fix $\varepsilon > 0$ and, by Lemma 6.2 (57), choose $N$ such that $\mathbb{P}^\gamma(\gamma \tau_N \leq T) \leq \varepsilon/2$ for any $\gamma \leq 1$. For any $n \leq N$,

$$\mathbb{P}^\gamma(\rho_{n+1} - \rho_n < \delta/\gamma) \leq \mathbb{P}^\gamma(\rho_n - \tau_n > \delta/\gamma) + \mathbb{P}^\gamma(\tau_{n+1} - \tau_n < 2\delta/\gamma).$$

By Lemma 6.2 (56) and Lemma 6.3(b), one can choose $\delta$ such that the second term is bounded by $\varepsilon/2N$ uniformly in $\gamma$. Then, by Lemma 6.4(a) and Lemma 6.3(b), the first term goes to 0 when $\gamma \to 0$. Therefore, for any $\varepsilon > 0$, there exists $\delta > 0$ such that

$$\limsup_{\gamma \to 0} \mathbb{P}^\gamma(\exists n \geq 0 : \rho_{n+1} - \rho_n < \delta/\gamma \text{ and } \rho_{n+1} \leq T/\gamma) \leq \varepsilon,$$

which implies (73). □

6.3. *Proof of Theorem* 4.1. Since $v_{n,m} \equiv 0$ on $\Omega_1 \cup \Omega_2$, we will always assume $n, m \geq 1$. Recall the definition (10) of the transition probabilities $\pi_{ij}(x,y)$ ($i,j \in \mathbb{N}^2$) of the embedded Markov chain $(\hat{X}_n, \hat{Y}_n; n \geq 0)$ associated to the $(x,y)$-type $GL$-population $(X_t, Y_t; t \geq 0)$ without mutation. Recall also the notation $\pi_{i_1,\ldots,i_k}(x,y)$ for the product $\pi_{i_1 i_2}(x,y)\pi_{i_2 i_3}(x,y)\cdots\pi_{i_{k-1}i_k}(x,y)$ and $S_{(n,m)\to \Gamma}$ for the set of paths linking $(n,m)$ to a subset $\Gamma$ of $\mathbb{N}^2$ without hitting $\Omega_1 \cup \Omega_2$ before $\Gamma$, that is, the set of all $k$-tuples $(i_1, i_2, \ldots, i_k)$ for all $k \geq 1$ such that $i_1 = (n,m)$, $i_2, \ldots, i_{k-1} \in \mathbb{N}^2 \setminus (\Omega_1 \cup \Omega_2)$ and $i_k \in \Gamma$.

Now,
$$u_{n,m}(x,y) = \sum_{k\geq 2, (i_1,\ldots,i_k)\in S_{(n,m)\to\Omega_1}} \pi_{i_1,\ldots,i_k}(x,y),$$

so if we prove that
$$R_{n,m}(x,y) := \sum_{k\geq 2, (i_1,\ldots,i_k)\in S_{(n,m)\to\Omega_1}} \left|\frac{\partial \pi_{i_1,\ldots,i_k}}{\partial y}(x,y)\right|$$

is finite, we get the differentiability of $u_{n,m}(x,y)$ and the inequality $|\partial u_{n,m}(x,y)/\partial y| \leq R_{n,m}(x,y)$. Observe that

$$R_{n,m}(x,y) \leq C_1 \sum_{k\geq 2, (i_1,\ldots,i_k)\in S_{(n,m)\to\Omega_1}} \sum_{l=1}^{k-1} \pi_{i_1,\ldots,i_l}(x,y)\pi_{i_{l+1},\ldots,i_k}(x,y),$$

where $C_1$ is defined in assumption (11) and $\pi_{i,i}(x,y) = 1$ by convention. Next, with $\|\cdot\|$ denoting the $L^1$-norm in $\mathbb{Z}^2$,

$$R_{n,m}(x,y) \leq C_1 \sum_{l\geq 1} \sum_{(n',m')\in(\mathbb{N}^*)^2} \sum_{(i_1,\ldots,i_l)\in S_{(n,m)\to(n',m')}} \pi_{i_1,\ldots,i_l}(x,y)$$
$$\times \sum_{\|\varepsilon\|=1} \sum_{k'\geq 0, (j_1,\ldots,j_{k'})\in S_{(n',m')+\varepsilon\to\Omega_1}} \pi_{j_1,\ldots,j_{k'}}(x,y),$$

with the convention that $k' = 0$ if $(n',m') + \varepsilon \in \Omega_1$, so that, if $\hat{T}_\Gamma$ denotes the first hitting of $\Gamma \subset \mathbb{N}^2$ by $(\hat{X}, \hat{Y})$,

$$R_{n,m}(x,y) \leq C_1 \sum_{l\geq 1} \sum_{(n',m')\in(\mathbb{N}^*)^2} \sum_{(i_1,\ldots,i_l)\in S_{(n,m)\to(n',m')}} \pi_{i_1,\ldots,i_l}(x,y)$$
$$\times \sum_{\|\varepsilon\|=1} \mathbb{P}_{(n',m')+\varepsilon}(\hat{T}_{\Omega_1} < \hat{T}_{\Omega_2})$$
$$\leq 4C_1 \sum_{l\geq 1} \mathbb{P}_{n,m}(\hat{T} > l)$$
$$= 4C_1 \mathbb{E}_{n,m}(\hat{T} - 1).$$

The proof of (12) is completed thanks to Proposition 6.1(a).

(b) Following exactly the same method as above, we obtain from (54) that

$$\sum_{(i_1,\ldots,i_k)\in S_{(n,m)\to\Omega_1}} \left|\frac{\partial^2 \pi_{i_1,\ldots,i_k}}{\partial y^2}(x,y)\right|$$

$$\leq 4C_1^2 \sum_{l\geq 1} \sum_{(n',m')\in(\mathbb{N}^*)^2} \sum_{(i_1,\ldots,i_l)\in S_{(n,m)\to(n',m')}} \pi_{i_1,\ldots,i_l}(x,y)$$

$$\times \sum_{\|\varepsilon\|=1} \mathbb{E}_{(n',m')+\varepsilon}(\hat{T}-1)$$

$$+ 4C_2 \mathbb{E}_{n,m}(\hat{T}-1)$$

$$\leq 16C_1^2 C \sum_{l\geq 1} \sum_{(n',m')\in(\mathbb{N}^*)^2} \sum_{(i_1,\ldots,i_l)\in S_{(n,m)\to(n',m')}} \pi_{i_1,\ldots,i_l}(x,y)(n'+m')$$

$$+ 4C_2 C(n+m)$$

$$\leq 16C_1^2 C \sum_{l\geq 1}(n+m+l)\mathbb{P}_{n,m}(\hat{T}>l) + 4C_2 C(n+m)$$

$$\leq 16C_1^2 C\mathbb{E}_{n,m}((\hat{T}-1)(n+m+\hat{T}/2)) + 4C_2 C(n+m),$$

and the result follows again from Proposition 6.1(a).

6.4. *Proof of Theorem* 4.2. We will use the classical method of tightness and martingale problem formulation to prove this theorem (e.g., [18]). We divide the proof into three steps:

*Step* 1. *Uniqueness of the limit process.* Strong existence and uniqueness for the SDE (15) follow standardly from the Lipschitz-continuity of its coefficients.

*Step* 2. *Tightness of the family of laws of* $Z^\epsilon$. For any $\epsilon > 0$, let $N_\epsilon(dh, du, dt)$ be a Poisson point process on $\mathbb{R}^k \times [0,1] \times \mathbb{R}_+$ with intensity measure $q_\epsilon(dh, du, dt) = \bar{M}(h)\,dh\bar{\beta}\bar{\chi}\,du\,dt/\epsilon^2$, where $\bar{\beta}$ and $\bar{\chi}$ are constants bounding the functions $\beta$ and $\chi$ from above, respectively, and $\bar{M}$ has been defined as the integrable function bounding the density $m(x,\cdot)$ of $M(x,\cdot)$ for any $x \in \mathbb{R}^k$. Then it is straightforward to check that $A_\epsilon$ is the infinitesimal generator of the Markov process $Z^\epsilon$:

$$Z_t^\epsilon = Z_0^\epsilon + \epsilon \int_0^t \int_0^1 \int_{\mathbb{R}^k} h\mathbf{1}_{\{u\leq (\beta(Z_{s-}^\epsilon)/\bar{\beta})(\chi(Z_{s-}^\epsilon,Z_{s-}^\epsilon+\epsilon h)/\bar{\chi})(m(Z_{s-}^\epsilon,h)/\bar{M}(h))\}}$$

$$\times N_\epsilon(dh, du, ds).$$

Since $\beta$ and $\chi$ are bounded, a process generated by $A_\epsilon$ is unique in law (e.g., [10]), and this construction characterizes the law of the process $Z^\epsilon$ appearing in the statement of Theorem 4.2. Let us denote this law by $\mathbf{P}_\epsilon$.

Observe that, if we denote by $\tilde{N}_\epsilon$ the compensated Poisson measure $N_\epsilon - q_\epsilon$, $Z_t^\epsilon$ can be decomposed as $Z_0^\epsilon + \tilde{Z}_t^\epsilon + \hat{Z}_t^\epsilon$, where

$$\tilde{Z}_t^\epsilon = \epsilon \int_0^t \int_0^1 \int_{\mathbb{R}^k} h \mathbf{1}_{\{u \leq (\beta(Z_{s-}^\epsilon)/\bar{\beta})(\chi(Z_{s-}^\epsilon, Z_{s-}^\epsilon + \epsilon h)/\bar{\chi})(m(Z_{s-}^\epsilon, h)/\bar{M}(h))\}} \tilde{N}_\epsilon(dh, du, ds)$$

and

$$\hat{Z}_t^\epsilon = \frac{1}{\epsilon} \int_0^t \int_{\mathbb{R}^k} h \beta(Z_{s-}^\epsilon) \chi(Z_{s-}^\epsilon, Z_{s-}^\epsilon + \epsilon h) M(Z_{s-}^\epsilon, dh) \, ds$$
$$= \frac{1}{\epsilon} \int_0^t \int_{\mathbb{R}^k} h \beta(Z_{s-}^\epsilon) [\chi(Z_{s-}^\epsilon, Z_{s-}^\epsilon + \epsilon h) - \chi(Z_{s-}^\epsilon, Z_{s-}^\epsilon)] M(Z_{s-}^\epsilon, dh) \, ds,$$

where the last equality follows from the fact that the mutation step law $M(x, \cdot)$ has 0 expectation.

We will use Aldous' criterion [1] to prove the tightness of the family of probability measures $(\mathbf{P}_\epsilon)_{\epsilon > 0}$ on $\mathbb{D}(\mathbb{R}_+, \mathbb{R}^k)$. Fix $\delta, \epsilon > 0$ and let $\tau$ and $\tau'$ be two stopping times such that $\tau < \tau' < \tau + \delta$. Since $|\chi(x, x + \epsilon h) - \chi(x, x)| \leq \epsilon K \|h\|$ for some constant $K$, $\|\hat{Z}_{\tau'}^\epsilon - \hat{Z}_\tau^\epsilon\| \leq \delta \bar{\beta} K M_2$, where $M_2 = \sup_x \int \|h\|^2 M(x, dh)$, which is finite by assumption. By standard results on stochastic integrals with respect to Poisson point measures,

$$\mathbf{E}_\epsilon(\|\tilde{Z}_{\tau'}^\epsilon - \tilde{Z}_\tau^\epsilon\|^2)$$
$$= \mathbf{E}_\epsilon\left(\int_\tau^{\tau'} \int_0^1 \int_{\mathbb{R}^k} \epsilon^2 \|h\|^2 \mathbf{1}_{\{u \leq (\beta(Z_{s-}^\epsilon)/\bar{\beta})(\chi(Z_{s-}^\epsilon, Z_{s-}^\epsilon + \epsilon h)/\bar{\chi})(m(Z_{s-}^\epsilon, h)/\bar{M}(h))\}} \right.$$
$$\left. \times q_\epsilon(dh, du, ds)\right)$$
$$\leq \delta \bar{\beta} \bar{\chi} M_2.$$

Therefore, for any $\eta > 0$,

$$\mathbf{P}_\epsilon(\|Z_{\tau'}^\epsilon - Z_\tau^\epsilon\| \geq \eta) \leq \mathbf{P}_\epsilon\left(\|\hat{Z}_{\tau'}^\epsilon - \hat{Z}_\tau^\epsilon\| \geq \frac{\eta}{2}\right) + \mathbf{P}_\epsilon\left(\|\tilde{Z}_{\tau'}^\epsilon - \tilde{Z}_\tau^\epsilon\| \geq \frac{\eta}{2}\right)$$
$$\leq \mathbf{1}_{\{2\delta\bar{\beta}KM_2 \geq \eta\}} + \frac{4\delta\bar{\beta}\bar{\chi}M_2}{\eta^2},$$

which converges to 0 as $\delta \to 0$, uniformly w.r.t. the choice of $\tau$ and $\tau'$. This is the first part of Aldous' criterion. For the second part, we have to prove the tightness of $(Z_t^\epsilon)_{\epsilon > 0}$ for any $t \geq 0$. Similar computations as above prove easily that $(\tilde{Z}_t^\epsilon)_{\epsilon > 0}$ and $(\hat{Z}_t^\epsilon)_{\epsilon > 0}$ are tight, and the tightness of $(Z_0^\epsilon)_{\epsilon > 0}$ follows from the assumption that it is bounded in $L^1$.

*Step* 3. *Martingale problem*. Let $\mathbf{P}_0$ be an accumulation point of $(\mathbf{P}_\epsilon)_{\epsilon > 0}$ as $\epsilon \to 0$ on $\mathbb{D}(\mathbb{R}_+, \mathbb{R}^k)$, endowed with the canonical filtration $\mathcal{F}_t$. Since the martingale problem for (15) is well posed, it suffices to show that, for any $\varphi \in \mathcal{C}^2(\mathcal{X})$, under $\mathbf{P}_0$, the process

$$M_t^\varphi(w) = \varphi(w_t) - \varphi(w_0) - \int_0^t A_0 \varphi(w_s) \, ds$$

on $\mathbb{D}(\mathbb{R}_+, \mathcal{X})$ is a local $\mathcal{F}_t$-martingale. We already know that, under $\mathbf{P}_\epsilon$,

$$M_t^{\epsilon,\varphi}(w) = \varphi(w_t) - \varphi(w_0) - \int_0^t A_\epsilon \varphi(w_s) \, ds$$

is a local $\mathcal{F}_t$-martingale. Since $\beta$ and $\chi$ are bounded, this is a square-integrable martingale as soon as $\varphi \in \mathcal{C}_b^3$.

Fix $\varphi \in \mathcal{C}_b^3$, $s > 0$ and $t > s$, and consider $p$ real numbers $0 \leq t_1 < \cdots < t_p \leq s$ for some $p \geq 1$, and a continuous bounded function $q : (\mathbb{R}^k)^p \to \mathbb{R}$. We can write

$$\left| \mathbf{E}_0 \left\{ q(w_{t_1}, \ldots, w_{t_p}) \left[ \varphi(w_t) - \varphi(w_s) - \int_s^t A_0 \varphi(w_u) \, du \right] \right\} \right|$$
$$\leq \left| \mathbf{E}_\epsilon \left\{ q(w_{t_1}, \ldots, w_{t_p}) \left[ \varphi(w_t) - \varphi(w_s) - \int_s^t A_\epsilon \varphi(w_u) \, du \right] \right\} \right|$$
$$+ \left| \mathbf{E}_\epsilon \left\{ q(w_{t_1}, \ldots, w_{t_p}) \int_s^t (A_\epsilon \varphi(w_u) - A_0 \varphi(w_u)) \, du \right\} \right|$$
$$+ \left| \mathbf{E}_0 \left\{ q(w_{t_1}, \ldots, w_{t_p}) \left[ \varphi(w_t) - \varphi(w_s) - \int_s^t A_0 \varphi(w_u) \, du \right] \right\} \right.$$
$$\left. - \mathbf{E}_\epsilon \left\{ q(w_{t_1}, \ldots, w_{t_p}) \left[ \varphi(w_t) - \varphi(w_s) - \int_s^t A_0 \varphi(w_u) \, du \right] \right\} \right|.$$

The first term of the r.h.s. is 0 since $M^{\epsilon,\varphi}$ is a $\mathbf{P}_\epsilon$-martingale. Because of the uniform convergence of generators (14), the second term converges to 0 as $\epsilon \to 0$. The third term also goes to 0 as $\epsilon \to 0$ because $\mathbf{P}_\epsilon$ converges to $\mathbf{P}_0$. Finally, since the l.h.s. does not depend on $\epsilon$, it is 0.

A classical use of the monotone class theorem allows to extend this equality to all $\mathcal{F}_s$-measurable bounded functions $q$, so $M^\varphi$ is a $\mathbf{P}_0$-martingale. This result can easily be extended to any $\mathcal{C}^2$ function $\varphi$ by a standard truncation technique, which completes the proof of Theorem 4.2.

Weierstrass Institute for Applied
  Analysis and Stochastics
Mohrenstrasse 39
10117 Berlin
Germany
E-mail: champagn@wias-berlin.de

Unit of Mathematical Evolutionary Biology
UMR 7625 Laboratoire d'Écologie
École Normale Supérieure
46 rue d'Ulm
F-75230 Paris Cedex 05
and
UMR 7625 Laboratoire d'Écologie
Université Pierre et Marie Curie-Paris 6
7 quai Saint Bernard
F-75252 Paris Cedex 05
France
E-mail: amaury.lambert@ens.fr